\documentclass[12pt]{article}

\usepackage{amsmath,amssymb,amsthm,amscd,a4wide}

\begin{document}

\newcommand{\End}{{\rm{End}\ts}}
\newcommand{\Hom}{{\rm{Hom}}}
\newcommand{\Mat}{{\rm{Mat}}}
\newcommand{\ch}{{\rm{ch}\ts}}
\newcommand{\chara}{{\rm{char}\ts}}
\newcommand{\diag}{ {\rm diag}}
\newcommand{\non}{\nonumber}
\newcommand{\wt}{\widetilde}
\newcommand{\wh}{\widehat}
\newcommand{\ot}{\otimes}
\newcommand{\la}{\lambda}
\newcommand{\La}{\Lambda}
\newcommand{\De}{\Delta}
\newcommand{\al}{\alpha}
\newcommand{\be}{\beta}
\newcommand{\ga}{\gamma}
\newcommand{\Ga}{\Gamma}
\newcommand{\ep}{\epsilon}
\newcommand{\ka}{\kappa}
\newcommand{\vk}{\varkappa}
\newcommand{\si}{\sigma}
\newcommand{\vp}{\varphi}
\newcommand{\de}{\delta}
\newcommand{\ze}{\zeta}
\newcommand{\om}{\omega}
\newcommand{\ee}{\epsilon^{}}
\newcommand{\su}{s^{}}
\newcommand{\hra}{\hookrightarrow}
\newcommand{\ve}{\varepsilon}
\newcommand{\ts}{\,}
\newcommand{\vac}{\mathbf{1}}
\newcommand{\vacr}{|\tss 0\rangle}
\newcommand{\vacl}{\langle 0\tss |}
\newcommand{\di}{\partial}
\newcommand{\qin}{q^{-1}}
\newcommand{\tss}{\hspace{1pt}}
\newcommand{\Sr}{ {\rm S}}
\newcommand{\U}{ {\rm U}}
\newcommand{\BL}{ {\overline L}}
\newcommand{\BE}{ {\overline E}}
\newcommand{\BP}{ {\overline P}}
\newcommand{\AAb}{\mathbb{A}\tss}
\newcommand{\CC}{\mathbb{C}\tss}
\newcommand{\KK}{\mathbb{K}\tss}
\newcommand{\QQ}{\mathbb{Q}\tss}
\newcommand{\SSb}{\mathbb{S}\tss}
\newcommand{\ZZ}{\mathbb{Z}\tss}
\newcommand{\X}{ {\rm X}}
\newcommand{\Y}{ {\rm Y}}
\newcommand{\Z}{{\rm Z}}
\newcommand{\Ac}{\mathcal{A}}
\newcommand{\achi}{\Ac_{\chi}}
\newcommand{\bachi}{\overline\Ac_{\chi}}
\newcommand{\Lc}{\mathcal{L}}
\newcommand{\Mc}{\mathcal{M}}
\newcommand{\Pc}{\mathcal{P}}
\newcommand{\Qc}{\mathcal{Q}}
\newcommand{\Tc}{\mathcal{T}}
\newcommand{\Sc}{\mathcal{S}}
\newcommand{\Bc}{\mathcal{B}}
\newcommand{\Dc}{\mathcal{D}}
\newcommand{\Ec}{\mathcal{E}}
\newcommand{\Fc}{\mathcal{F}}
\newcommand{\Hc}{\mathcal{H}}
\newcommand{\Uc}{\mathcal{U}}
\newcommand{\Vc}{\mathcal{V}}
\newcommand{\Wc}{\mathcal{W}}
\newcommand{\Yc}{\mathcal{Y}}
\newcommand{\Ar}{{\rm A}}
\newcommand{\Br}{{\rm B}}
\newcommand{\Ir}{{\rm I}}
\newcommand{\Fr}{{\rm F}}
\newcommand{\Jr}{{\rm J}}
\newcommand{\Or}{{\rm O}}
\newcommand{\GL}{{\rm GL}}
\newcommand{\Spr}{{\rm Sp}}
\newcommand{\Rr}{{\rm R}}
\newcommand{\Zr}{{\rm Z}}
\newcommand{\gl}{\mathfrak{gl}}
\newcommand{\middd}{{\rm mid}}
\newcommand{\ev}{{\rm ev}}
\newcommand{\Pf}{{\rm Pf}}
\newcommand{\Norm}{{\rm Norm\tss}}
\newcommand{\oa}{\mathfrak{o}}
\newcommand{\spa}{\mathfrak{sp}}
\newcommand{\osp}{\mathfrak{osp}}
\newcommand{\g}{\mathfrak{g}}
\newcommand{\h}{\mathfrak h}
\newcommand{\n}{\mathfrak n}
\newcommand{\z}{\mathfrak{z}}
\newcommand{\Zgot}{\mathfrak{Z}}
\newcommand{\p}{\mathfrak{p}}
\newcommand{\sll}{\mathfrak{sl}}
\newcommand{\agot}{\mathfrak{a}}
\newcommand{\qdet}{ {\rm qdet}\ts}
\newcommand{\Ber}{ {\rm Ber}\ts}
\newcommand{\HC}{ {\mathcal HC}}
\newcommand{\cdet}{ {\rm cdet}}
\newcommand{\tr}{ {\rm tr}}
\newcommand{\gr}{ {\rm gr}}
\newcommand{\str}{ {\rm str}}
\newcommand{\loc}{{\rm loc}}
\newcommand{\Gr}{{\rm G}}
\newcommand{\sgn}{ {\rm sgn}\ts}
\newcommand{\ba}{\bar{a}}
\newcommand{\bb}{\bar{b}}
\newcommand{\bi}{\bar{\imath}}
\newcommand{\bj}{\bar{\jmath}}
\newcommand{\bk}{\bar{k}}
\newcommand{\bl}{\bar{l}}
\newcommand{\hb}{\mathbf{h}}
\newcommand{\Sym}{\mathfrak S}
\newcommand{\fand}{\quad\text{and}\quad}
\newcommand{\Fand}{\qquad\text{and}\qquad}
\newcommand{\For}{\qquad\text{or}\qquad}
\newcommand{\OR}{\qquad\text{or}\qquad}

\renewcommand{\theequation}{\arabic{section}.\arabic{equation}}

\newtheorem{thm}{Theorem}[section]
\newtheorem{lem}[thm]{Lemma}
\newtheorem{prop}[thm]{Proposition}
\newtheorem{cor}[thm]{Corollary}
\newtheorem{conj}[thm]{Conjecture}
\newtheorem*{mthm}{Main Theorem}
\newtheorem*{mthma}{Theorem A}
\newtheorem*{mthmb}{Theorem B}

\theoremstyle{definition}
\newtheorem{defin}[thm]{Definition}

\theoremstyle{remark}
\newtheorem{remark}[thm]{Remark}
\newtheorem{example}[thm]{Example}

\newcommand{\bth}{\begin{thm}}
\renewcommand{\eth}{\end{thm}}
\newcommand{\bpr}{\begin{prop}}
\newcommand{\epr}{\end{prop}}
\newcommand{\ble}{\begin{lem}}
\newcommand{\ele}{\end{lem}}
\newcommand{\bco}{\begin{cor}}
\newcommand{\eco}{\end{cor}}
\newcommand{\bde}{\begin{defin}}
\newcommand{\ede}{\end{defin}}
\newcommand{\bex}{\begin{example}}
\newcommand{\eex}{\end{example}}
\newcommand{\bre}{\begin{remark}}
\newcommand{\ere}{\end{remark}}
\newcommand{\bcj}{\begin{conj}}
\newcommand{\ecj}{\end{conj}}

\newcommand{\bal}{\begin{aligned}}
\newcommand{\eal}{\end{aligned}}
\newcommand{\beq}{\begin{equation}}
\newcommand{\eeq}{\end{equation}}
\newcommand{\ben}{\begin{equation*}}
\newcommand{\een}{\end{equation*}}

\newcommand{\bpf}{\begin{proof}}
\newcommand{\epf}{\end{proof}}

\def\beql#1{\begin{equation}\label{#1}}

\title{\Large\bf Higher Sugawara operators for the quantum affine algebras of type $A$}

\author{{Luc Frappat,\ \   Naihuan Jing,\ \   Alexander Molev\ \   and\ \   Eric Ragoucy}}

\date{} 
\maketitle

\vspace{4 mm}

\begin{abstract}
We give explicit formulas for the elements of the center of
the completed quantum affine algebra in type $A$
at the critical level which are associated with the fundamental representations.
We calculate the images of these elements
under a Harish-Chandra-type homomorphism. These images
coincide with those in the free field realization
of the quantum affine algebra and
reproduce generators of the $q$-deformed classical
$\Wc$-algebra of Frenkel and Reshetikhin.

\medskip

Preprint LAPTH-027/15

%

\end{abstract}


\vspace{5 mm}

%


\section{Introduction}
\label{sec:int}
\setcounter{equation}{0}

Let $\g$ be a simple Lie algebra over $\CC$ and let $\U_q(\wh\g)$ denote the quantum
affine algebra associated with $\g$.
Due to the work of Reshetikhin and Semenov-Tian-Shansky~\cite{rs:ce},
to every finite-dimensional representation $V$
of $\U_q(\wh\g)$ one can associate a formal Laurent series $\ell_V(z)$ whose coefficients
belong to the center $\Zr_q(\wh\g)$ of the completion $\wt\U_q(\wh\g)_{\text{\rm cri}}$
of the quantum affine algebra at the critical level. The map $V\mapsto \ell_V(z)$ was further
studied by Ding and Etingof~\cite{de:cq} who showed that if the coefficients of
$\ell_V(z)$ are regarded as operators on highest weight modules at the critical level, then
it possesses properties of a homomorphism from
the Grothendieck ring $\text{\rm Rep}\ts \U_q(\wh\g)$ to formal series in $z$.
Furthermore, the coefficients of $\ell_V(z)$ were shown
to generate all singular vectors in Verma modules \cite{de:cq}.
This relied upon connections of the series $\ell_V(z)$ with transfer matrices; see also work
of Frenkel and Reshetikhin~\cite[Sec.~8]{fr:qc} for its relationship with
the $q$-characters of finite-dimensional representations of $\U_q(\wh\g)$.

By a conjecture of Frenkel and Reshetikhin~\cite{fr:qa}, \cite{fr:qc}, the center
$\Zr_q(\wh\g)$ is isomorphic to the $q$-deformed
classical $\Wc$-algebra, as a Poisson algebra. More precisely, Conjecture~1 of
\cite{fr:qa} applies to the completion of the
central subalgebra of $\wt\U_q(\wh\g)_{\text{\rm cri}}$
generated by the coefficients of the series $\ell_{V_1}(z),\dots,\ell_{V_n}(z)$ associated
with all fundamental representations $V_i$ of $\U_q(\wh\g)$. Its proof was sketched in \cite{fr:qa}
for $\g=\sll_N$. The isomorphism is provided by the free field
(or Wakimoto) realization of the quantum affine algebra due to
Awata, Odake and Shiraishi~\cite{aos:fb} which extended
an earlier work \cite{fj:vr} on the vertex representations
from level $1$ to an arbitrary level.

The results of \cite{fr:qa}
generalize the Feigin--Frenkel theorem which establishes a Poisson
algebra isomorphism between the center $\z(\wh\g)$ of the affine
vertex algebra $V(\g)$ at the critical level
and the classical $\Wc$-algebra associated with the Langlands dual Lie algebra ${}^L\g$;
see \cite{f:lc} for a detailed exposition.
The Feigin--Frenkel
center $\z(\wh\g)$ is an algebra of polynomials which can be identified
with a commutative subalgebra of the universal
enveloping algebra $\U\big(t^{-1}\g[t^{-1}]\big)$.
As discovered by Feigin, Frenkel and Reshetikhin~\cite{ffr:gm},
the higher degree Hamiltonians in the Gaudin model can be obtained from
generators of $\z(\wh\g)$; see also \cite{fft:gm} and \cite{r:si}.
Explicit constructions of generators of $\z(\wh\g)$
were given in \cite{cm:ho}, \cite{ct:qs}
and \cite{m:ff} for types $A, B, C$ and $D$; see also \cite{mm:yc}
for their images in the classical $\Wc$-algebras
and \cite{mr:mm} for super-analogues of these constructions with $\g=\gl_{m|n}$.

Our goal in this paper is to give similar explicit formulas for
higher degree Sugawara operators for $\U_q(\wh\gl_n)$; i.e., for elements
of the center of $\wt\U_q(\wh\gl_n)_{\text{\rm cri}}$ (Theorem~\ref{thm:centr}).
The formulas express the operators in terms of the $RLL$-presentation
$\U_q(\wh\gl_n)$; see \cite{fri:qa} and \cite{rs:ce}.
We use a version of the Poincar\'e--Birkhoff--Witt theorem
for this presentation to introduce an analogue of the Harish-Chandra homomorphism
and calculate the images of the central elements under the homomorphism.

Then we apply the Ding--Frenkel isomorphism~\cite{df:it}
between the $RLL$ and
Drinfeld presentations to calculate the images of the Sugawara operators
in the $q$-deformed classical $\Wc$-algebra by using the approach of \cite{fr:qa}
based on the free field realization \cite{aos:fb}.
Our generators correspond to the fundamental representations
of $\U_q(\wh\gl_n)$ and essentially coincide with those of \cite{fr:qa}
up to an appropriate identification of the parameters. The construction
involves a fusion formula for the $q$-deformed antisymmetrizer
expressing it in terms of the trigonometric $R$-matrices.

As an application of Theorem~\ref{thm:centr}, we produce explicit invariants
of the $q$-analogue $V_q(\gl_n)$ of the vacuum module over the quantum affine
algebra. The invariants are
obtained by the action of the higher degree Sugawara operators
on the vacuum vector; cf. \cite{fr:qc}.

\section{Quantum affine algebra}
\label{sec:qaa}
\setcounter{equation}{0}

We use the $RLL$ presentation of $\U_q(\wh\gl_n)$ introduced in \cite{rs:ce}; see also
\cite{fri:qa}.
We regard $q$ as a nonzero complex number which is not
a root of unity.
Introduce the two-parameter $R$-matrix $R(u,v)\in\End\CC^n\ot\End\CC^n$ by
\begin{align}
R(u,v)={}&(u-v)\sum_{i\ne j}e_{ii}\ot e_{jj}+(\qin u-q\tss v)
\sum_{i}e_{ii}\ot e_{ii}
\non\\
{}+ {}&(\qin-q)\tss u\tss\sum_{i> j}e_{ij}\ot
e_{ji}+ (\qin-q)\tss v\tss\sum_{i< j}e_{ij}\ot e_{ji},
\label{ruv}
\end{align}
where $e_{ij}\in \End\CC^n$ are the standard matrix units. We will also need
the one-parameter $R$-matrices
\begin{align}
\overline R(x)=\frac{R(x,1)}{\qin x-q}={}&\sum_{i}e_{ii}\ot e_{ii}
+\frac{1-x}{q-\qin x}\ts\sum_{i\ne j}e_{ii}\ot e_{jj}
\non\\[0.4em]
{}+{}&\frac{(q-\qin)\ts x}{q-\qin x}\ts\sum_{i> j}e_{ij}\ot
e_{ji}+ \frac{q-\qin}{q-\qin x}\ts\sum_{i< j}e_{ij}\ot e_{ji},
\label{rx}
\end{align}
and
\beql{rf}
R(x)=f(x)\tss \overline R(x),
\eeq
where
\ben
f(x)=1+\sum_{k=1}^{\infty}f_kx^k,\qquad f_k=f_k(q),
\een
is a formal power series in $x$ whose coefficients $f_k$ are
rational functions in $q$ uniquely determined
by the relation
\beql{fx}
f(xq^{2n})=f(x)\ts\frac{(1-xq^2)\tss(1-xq^{2n-2})}{(1-x)\tss(1-xq^{2n})}.
\eeq
They can be found by the recurrence
\ben
f_k=-\frac{(1-q^2)(1-q^{2n-2})}{1-q^{2n}}
\sum_{i=1}^{k}\frac{1-q^{2ni}}{1-q^{2nk}}\ts f_{k-i},\qquad k\geqslant 1,
\een
with $f_0=1$.
Equivalently, $f(x)$ can be given by
\ben
f(x)=\frac{(x;q^{2n})_{\infty}\tss(xq^{2n};q^{2n})_{\infty}}
{(xq^{2};q^{2n})_{\infty}\tss(xq^{2n-2};q^{2n})_{\infty}},
\qquad (a;b)_{\infty}:=\prod_{r=0}^{\infty}(1-a\tss b^r),
\een
where the coefficients of the powers $x^k$ are power series in $q$ converging to $f_k$ for $|q|<1$.

The {\em quantum affine algebra $\U_q(\wh\gl_n)$}
is generated by elements
\ben
l^+_{ij}[-r],\qquad l^-_{ij}[r]\qquad\text{with}\quad 1\leqslant i,j\leqslant n,\qquad r=0,1,\dots,
\een
and the invertible central element $q^c$,
subject to the defining relations
\begin{align}
l^+_{ji}[0]&=l^-_{ij}[0]=0\qquad&&\text{for}\qquad 1\leqslant i<j\leqslant n,
\label{llt}\\
l^+_{ii}[0]\ts l^-_{ii}[0]&=l^-_{ii}[0]\ts l^+_{ii}[0]=1\qquad&&\text{for}\qquad i=1,\dots,n,
\label{lplm}
\end{align}
and
\begin{align}
R(u/v)L_1^{\pm}(u)L_2^{\pm}(v)&=L_2^{\pm}(v)L_1^{\pm}(u)R(u/v),
\label{RLL}\\[0.2em]
R(uq^{-c}/v)L_1^{+}(u)L_2^{-}(v)&=L_2^{-}(v)L_1^{+}(u)R(uq^{c}/v).
\label{RLLpm}
\end{align}
In the last two relations we consider the matrices $L^{\pm}(u)=\big[\tss l^{\pm}_{ij}(u)\big]$,
whose entries are formal power series in $u$ and $u^{-1}$,
\beql{serlpm}
l^{+}_{ij}(u)=\sum_{r=0}^{\infty}l^{+}_{ij}[-r]\tss u^r,\qquad
l^{-}_{ij}(u)=\sum_{r=0}^{\infty}l^{-}_{ij}[r]\tss u^{-r}.
\eeq
Here and below we regard the matrices as elements
\ben
L^{\pm}(u)=\sum_{i,j=1}^n e_{ij}\ot l^{\pm}_{ij}(u)\in\End\CC^n\ot\U_q(\wh\gl_n)[[u^{\pm1}]]
\een
and use a subscript to indicate a copy of the matrix in the multiple
tensor product algebra
\beql{multtpr}
\underbrace{\End\CC^n\ot\dots\ot\End\CC^n}_k\ot\U_q(\wh\gl_n)[[u^{\pm1}]]
\eeq
so that
\beql{la}
L^{\pm}_a(u)=\sum_{i,j=1}^n 1^{\ot (a-1)}\ot e_{ij}\ot 1^{\ot (k-a)}\ot l^{\pm}_{ij}(u).
\eeq
In particular, we take $k=2$ for the defining relations \eqref{RLL} and \eqref{RLLpm}.

This notation for elements of algebras of the form \eqref{multtpr} will be extended
as follows. For an element
\ben
C=\sum_{i,j,r,s=1}^n c^{}_{ijrs}\ts e_{ij}\ot e_{rs}\in
\End \CC^n\ot\End \CC^n,
\een
and any two indices $a,b\in\{1,\dots,k\}$ such that $a\ne b$,
we denote by $C_{a\tss b}$ the element of the algebra $(\End\CC^n)^{\ot k}$ with $k\geqslant 2$
given by
\beql{cab}
C_{a\tss b}=\sum_{i,j,r,s=1}^n c^{}_{ijrs}\ts (e_{ij})_a\tss (e_{rs})_b,
\qquad
(e_{ij})_a=1^{\ot(a-1)}\ot e_{ij}\ot 1^{\ot(k-a)}.
\eeq
We regard the matrix transposition as the linear map
\ben
t:\End\CC^n\to\End\CC^n,\qquad e_{ij}\mapsto e_{ji}.
\een
For any $a\in\{1,\dots,k\}$ we will denote by $t_a$ the corresponding
partial transposition on the algebra \eqref{multtpr} which acts as $t$ on the
$a$-th copy of $\End \CC^n$ and as the identity map on all the other tensor factors.

The $R$-matrix \eqref{rf} satisfies the following {\em crossing symmetry relations} \cite{fri:qa}:
\beql{cs}
\big(R_{12}(x)^{-1}\big)^{t_2} D_2 R_{12}(xq^{2n})^{t_2}=D_2
\Fand
R_{12}(xq^{2n})^{t_1}\tss D_1\big(R_{12}(x)^{-1}\big)^{t_1}=D_1,
\eeq
where $D$ denotes the diagonal $n\times n$ matrix
\beql{d}
D=\diag\big[q^{n-1}, q^{n-3},\dots, q^{-n+1}\big]
\eeq
with the meaning of subscripts as in \eqref{cab}.

\section{Main theorem}
\label{sec:mt}
\setcounter{equation}{0}

Denote by $\U_q(\wh\gl_n)_{\text{\rm cri}}$ the quantum affine algebra
{\em at the critical level $c=-n$}, which is
the quotient of $\U_q(\wh\gl_n)$ by the relation $q^c=q^{-n}$.
Its completion $\wt\U_q(\wh\gl_n)_{\text{\rm cri}}$ is defined as
the inverse limit
\beql{compl}
\wt\U_q(\wh\gl_n)_{\text{\rm cri}}=\lim_{\longleftarrow}
\U_q(\wh\gl_n)_{\text{\rm cri}}/J_p, \qquad p>0,
\eeq
where $J_p$ denotes the left ideal of $\U_q(\wh\gl_n)_{\text{\rm cri}}$ generated by all elements
$l^{-}_{ij}[r]$ with $r\geqslant p$.
Elements of the center $\Zr_q(\wh\gl_n)$ of $\wt\U_q(\wh\gl_n)_{\text{\rm cri}}$
are known as {\em Sugawara operators}.

Consider the $q$-permutation operator
$P^{\tss q}\in\End(\CC^n\ot\CC^n)\cong\End\CC^n\ot\End\CC^n$
defined by
\beql{qperm}
P^{\tss q}=\sum_{i}e_{ii}\ot e_{ii}+ q\tss\sum_{i> j}e_{ij}\ot
e_{ji}+ \qin\sum_{i< j}e_{ij}\ot e_{ji}.
\eeq
The symmetric group $\Sym_k$ acts on the space $(\CC^n)^{\ot\tss k}$
by $s_i\mapsto P^{\tss q}_{s_i}:=
P^{\tss q}_{i,i+1}$ for $i=1,\dots,k-1$,
where $s_i$ denotes the transposition $(i,i+1)$.
If $\si=s_{i_1}\cdots s_{i_l}$ is a reduced decomposition
of an element $\si\in \Sym_k$ we set
$P^{\tss q}_{\si}=P^{\tss q}_{s_{i_1}}\cdots P^{\tss q}_{s_{i_l}}$.
We denote by $A^{(k)}$ the image of the normalized antisymmetrizer
associated with the $q$-permutations:
\beql{antisym}
A^{(k)}=\frac{1}{k!}\ts\sum_{\si\in\Sym_k}\sgn\tss\si\cdot P^{\tss q}_{\si}
\eeq
so that
$\big(A^{(k)}\big)^2=A^{(k)}$.

For each $k=1,\dots,n$ introduce the Laurent series $\ell_k(z)$ in $z$ by
\beql{lkz}
\ell_k(z)=\tr^{}_{1,\dots,k}\ts A^{(k)}L^+_1(z)\dots L^+_k(zq^{-2k+2})
L^-_k(zq^{-n-2k+2})^{-1}\dots L^-_1(zq^{-n})^{-1} D_1\dots D_k,
\eeq
where $D$ is the diagonal matrix \eqref{d} and
the trace is taken over all $k$ copies of $\End\CC^n$ in \eqref{multtpr}.
All coefficients of the series $\ell_k(z)$
are elements of the algebra $\wt\U_q(\wh\gl_n)_{\text{\rm cri}}$.
We will also need an equivalent formula for $\ell_k(z)$ which is obtained by using
the following well-known particular case
of the fusion procedure for the $R$-matrix \eqref{ruv}; see \cite{c:ni}.

\ble\label{lem:deca}
Set $v_a=zq^{-2a+2}$ for $a=1,\dots,k$.
We have
\ben
\prod_{1\leqslant a<b\leqslant k}\ts R_{a b}(v_a,v_b)=k!\ts z^{k(k-1)/2}\ts
\prod_{0\leqslant a<b\leqslant k-1}
(q^{-2a}-q^{-2b})\ts A^{(k)},
\een
where the product is taken in the lexicographical order on the pairs $(a,b)$.
\qed
\ele

Applying \eqref{RLL} and Lemma~\ref{lem:deca} we obtain the relations
\begin{align}
A^{(k)}L^+_1(v_1)\dots L^+_k(v_k)
{}&=L^+_k(v_k)\dots L^+_1(v_1)\tss A^{(k)},
\label{lpa}\\[0.4em]
A^{(k)}L^-_k(v_kq^{-n})^{-1}\dots L^-_1(v_1q^{-n})^{-1}
{}&=L^-_1(v_1q^{-n})^{-1}\dots L^-_k(v_kq^{-n})^{-1}\tss A^{(k)}.
\label{lina}
\end{align}
We also have
\ben
R(u,v) D_1D_2=D_2D_1R(u,v)
\een
and $D_1D_2=D_2D_1$ so that
\beql{da}
A^{(k)}D_1\dots D_k=D_1\dots D_k\tss A^{(k)}.
\eeq
Thus, $\ell_k(z)$ defined in \eqref{lkz} can also be given by the formula
\beql{lkzinv}
\ell_k(z)=\tr^{}_{1,\dots,k}\ts L^+_k(v_k)\dots L^+_1(v_1)\tss
L^-_1(v_1q^{-n})^{-1}\dots L^-_k(v_kq^{-n})^{-1}D_1\dots D_k\tss A^{(k)}.
\eeq

The following is our main result which provides explicit formulas
for higher Sugawara operators.

\bth\label{thm:centr}
The coefficients of $\ell_k(z)$ belong to the center of the completed quantum
affine algebra at the critical level $\wt\U_q(\wh\gl_n)_{\text{\rm cri}}$ for all $k=1,\dots,n$.
\eth

\bpf
Introduce an extra copy of the endomorphism algebra $\End\CC^n$ in \eqref{multtpr}
and label it by $0$ to work with the algebra
\beql{multtprze}
\End\CC^n\ot\big(\End\CC^n\big)^{\ot k}\ot\wt\U_q(\wh\gl_n)_{\text{\rm cri}}.
\eeq
It will be sufficient to verify what $\ell_k(z)$ commutes with $L^{\pm}_0(u)$.
By using \eqref{RLL}
we get
\ben
\bal
L^+_0(u)&L^+_k(v_k)\dots L^+_1(v_1)\\[0.2em]
{}&=R_{0k}(u/v_k)^{-1}\dots R_{01}(u/v_1)^{-1}
R_{01}(u/v_1)\dots R_{0k}(u/v_k)L^+_0(u)L^+_k(v_k)\dots L^+_1(v_1)\\[0.2em]
{}&=R_{0k}(u/v_k)^{-1}\dots R_{01}(u/v_1)^{-1}
L^+_k(v_k)\dots L^+_1(v_1)L^+_0(u)R_{01}(u/v_1)\dots R_{0k}(u/v_k).
\eal
\een
Relation \eqref{RLLpm} implies
\ben
L^+_0(u)R_{0a}(u/v_a)L^-_a(v_aq^{-n})^{-1}=L^-_a(v_aq^{-n})^{-1}R_{0a}(u\tss q^{2n}/v_a)L^+_0(u),
\qquad a=1,\dots,k,
\een
and so
\begin{multline}
L^+_0(u)R_{01}(u/v_1)\dots R_{0k}(u/v_k)
L^-_1(v_1q^{-n})^{-1}\dots L^-_k(v_kq^{-n})^{-1}\\[0.3em]
{}=
L^-_1(v_1q^{-n})^{-1}\dots L^-_k(v_kq^{-n})^{-1}
R_{01}(u\tss q^{2n}/v_1)\dots R_{0k}(u\tss q^{2n}/v_k)L^+_0(u).
\non
\end{multline}
Thus, to conclude that $L^+_0(u)\tss\ell_k(z)=\ell_k(z)\tss L^+_0(u)$ we need
to show that the trace
\begin{multline}
\tr^{}_{1,\dots,k}\ts R_{0k}(u/v_k)^{-1}\dots R_{01}(u/v_1)^{-1}
L^+_k(v_k)\dots L^+_1(v_1)\\[0.3em]
{}\times L^-_1(v_1q^{-n})^{-1}\dots L^-_k(v_kq^{-n})^{-1}
R_{01}(u\tss q^{2n}/v_1)\dots R_{0k}(u\tss q^{2n}/v_k)
D_1\dots D_k\tss A^{(k)}
\label{trace}
\end{multline}
equals $\ell_k(z)$.
The $R$-matrix $R(u,v)$ satisfies the Yang--Baxter equation
\ben
R_{12}(u,v)R_{13}(u,w)R_{23}(v,w)=R_{23}(v,w)R_{13}(u,w)R_{12}(u,v)
\een
which implies
\ben
R_{12}(u,v)\tss R_{13}(u/w)\tss R_{23}(v/w)=R_{23}(v/w)\tss R_{13}(u/w)\tss R_{12}(u,v).
\een
Therefore, Lemma~\ref{lem:deca} gives
\begin{align}
R_{01}(u\tss q^{2n}/v_1)\dots R_{0k}(u\tss q^{2n}/v_k)\tss A^{(k)}
&=A^{(k)}R_{0k}(u\tss q^{2n}/v_k)\dots R_{01}(u\tss q^{2n}/v_1),
\label{ra}\\[0.3em]
R_{0k}(u/v_k)^{-1}\dots R_{01}(u/v_1)^{-1}A^{(k)}&
=A^{(k)}R_{01}(u/v_1)^{-1}\dots R_{0k}(u/v_k)^{-1}.
\label{rainv}
\end{align}
Applying \eqref{lpa}, \eqref{lina}, \eqref{da} and \eqref{ra} to the expression
under the trace in \eqref{trace}, we will bring it to the form
\begin{multline}
\tr^{}_{1,\dots,k}\ts R_{0k}(u/v_k)^{-1}\dots R_{01}(u/v_1)^{-1} A^{(k)}
L^+_1(v_1)\dots L^+_k(v_k)\\[0.3em]
{}\times L^-_k(v_kq^{-n})^{-1}\dots L^-_1(v_1q^{-n})^{-1}
R_{0k}(u\tss q^{2n}/v_k)\dots R_{01}(u\tss q^{2n}/v_1)
D_1\dots D_k.
\label{traceam}
\end{multline}
Now write $A^{(k)}=(A^{(k)})^2$ and move one copy of $A^{(k)}$ to the left by using
\eqref{rainv}, and move the other copy back to its right-most position.
We get
\begin{multline}
\tr^{}_{1,\dots,k}\ts A^{(k)}R_{01}(u/v_1)^{-1}\dots R_{0k}(u/v_k)^{-1}
L^+_k(v_k)\dots L^+_1(v_1)\\[0.3em]
{}\times L^-_1(v_1q^{-n})^{-1}\dots L^-_k(v_kq^{-n})^{-1}
R_{01}(u\tss q^{2n}/v_1)\dots R_{0k}(u\tss q^{2n}/v_k)
D_1\dots D_k\tss A^{(k)}.
\label{traceala}
\end{multline}
Use the cyclic property of trace to move the left copy of $A^{(k)}$ to
the right-most position and replace $(A^{(k)})^2$ with $A^{(k)}$.
As a result of these transformations, the order
of the first $k$ factors in \eqref{trace} will be reversed,
while the rest of the expression remains unchanged.
Therefore, we can also write it in the form \eqref{traceam}
with the order
of the first $k$ factors reversed; that is, of the form
$
\tr^{}_{1,\dots,k}\ts XY
$
with
\ben
X=R_{01}(u/v_1)^{-1}\dots R_{0k}(u/v_k)^{-1}
A^{(k)}
L^+_1(v_1)\dots L^+_k(v_k)
L^-_k(v_kq^{-n})^{-1}\dots L^-_1(v_1q^{-n})^{-1}
\een
and
\ben
Y=R_{0k}(u\tss q^{2n}/v_k)\dots R_{01}(u\tss q^{2n}/v_1)
D_1\dots D_k.
\een
Now use the property
\ben
\tr^{}_{1,\dots,k}\ts XY=\tr^{}_{1,\dots,k}\ts X^{t_1\dots t_k}Y^{t_1\dots t_k}.
\een
We have
\ben
X^{t_1\dots t_k}=L^{t_1\dots t_k}
\big(R_{01}(u/v_1)^{-1}\big)^{t_1}\dots \big(R_{0k}(u/v_k)^{-1}\big)^{t_k},
\een
where we have set
\ben
L=A^{(k)}
L^+_1(v_1)\dots L^+_k(v_k)
L^-_k(v_kq^{-n})^{-1}\dots L^-_1(v_1q^{-n})^{-1}.
\een
Furthermore,
\ben
Y^{t_1\dots t_k}=D_1\dots D_k
R_{0k}(u\tss q^{2n}/v_k)^{t_k}\dots R_{01}(u\tss q^{2n}/v_1)^{t_1}.
\een
Hence, the first crossing symmetry relation in \eqref{cs} gives
\ben
\bal
\tr^{}_{1,\dots,k}\ts X^{t_1\dots t_k}Y^{t_1\dots t_k}
&=\tr^{}_{1,\dots,k}\ts L^{t_1\dots t_k}D_1\dots D_k\\[0.3em]
&{}=\tr^{}_{1,\dots,k}\ts A^{(k)}
L^+_1(v_1)\dots L^+_k(v_k)
L^-_k(v_kq^{-n})^{-1}\dots L^-_1(v_1q^{-n})^{-1}
D_1\dots D_k
\eal
\een
which coincides with $\ell_k(z)$ as defined in \eqref{lkz}
thus completing the proof for $L^+_0(u)$.

The argument showing
that $L^-_0(u)\tss\ell_k(z)=\ell_k(z)\tss L^-_0(u)$ is quite similar,
so we only briefly outline it. Using \eqref{RLLpm}
we get
\ben
\bal
L^-_0(u)&L^+_k(v_k)\dots L^+_1(v_1)=R_{k\tss 0}(v_kq^n/u)\dots R_{10}(v_1q^n/u)\\[0.3em]
&{}\times L^+_k(v_k)\dots L^+_1(v_1)
L^-_0(u)R_{10}(v_1q^{-n}/u)^{-1}\dots R_{k\tss 0}(v_kq^{-n}/u)^{-1}.
\eal
\een
Next, due to \eqref{RLL}, we have
\begin{multline}
L^-_0(u)R_{10}(v_1q^{-n}/u)^{-1}\dots R_{k\tss 0}(v_kq^{-n}/u)^{-1}
L^-_1(v_1q^{-n})^{-1}\dots L^-_k(v_kq^{-n})^{-1}\\[0.3em]
{}=
L^-_1(v_1q^{-n})^{-1}\dots L^-_k(v_kq^{-n})^{-1}
R_{10}(v_1q^{-n}/u)^{-1}\dots R_{k\tss 0}(v_kq^{-n}/u)^{-1}L^-_0(u).
\non
\end{multline}
The argument is completed by verifying that the trace
\begin{multline}
\tr^{}_{1,\dots,k}\ts R_{k\tss 0}(v_kq^n/u)\dots R_{10}(v_1q^n/u)
L^+_k(v_k)\dots L^+_1(v_1)L^-_1(v_1q^{-n})^{-1}\dots L^-_k(v_kq^{-n})^{-1}
\\[0.3em]
{}\times
R_{10}(v_1q^{-n}/u)^{-1}\dots R_{k\tss 0}(v_kq^{-n}/u)^{-1}
D_1\dots D_k\tss A^{(k)}
\non
\end{multline}
coincides with $\ell_k(z)$. This is done in the same way as for the expression \eqref{trace}
with the use of Lemma~\ref{lem:deca} and the second crossing symmetry relation in \eqref{cs}.
\epf

As a corollary of Theorem~\ref{thm:centr}, we obtain an explicit description
of invariants of the vacuum module over the quantum affine algebra;
cf. \cite{cm:ho}, \cite{ct:qs}. By definition, the {\em vacuum module at the critical level}
$V_q(\gl_n)$ is the quotient of $\U_q(\wh\gl_n)_{\text{\rm cri}}$ by the left ideal
generated by all elements $l^-_{ij}[r]$ with $r>0$ and by
the elements $l^-_{i\tss j}[0]-\de_{i\tss j}$ with $i\geqslant j$.
The module $V_q(\gl_n)$ is generated by the vector $\vac$
(the image of $1\in \U_q(\wh\gl_n)_{\text{\rm cri}}$ in the quotient)
such that
\ben
L^-(u)\tss\vac=I\ts\vac,
\een
where $I$ denotes the identity matrix. As a vector space,
$V_q(\gl_n)$ can be identified with the subalgebra
$\Y_q(\gl_n)$ of $\U_q(\wh\gl_n)_{\text{\rm cri}}$ generated by
the coefficients of all series $l^+_{ij}(u)$
subject to the additional relations $l^+_{ii}[0]=1$.
This relies on the Poincar\'e--Birkhoff--Witt theorem for the quantum affine
algebra; see, e.g., Sec.~\ref{sec:pbw} below.
The subspace of invariants of $V_q(\gl_n)$ is defined by
\ben
\z_q(\wh\gl_n)=\{v\in V_q(\gl_n)\ |\ L^-(u)\tss v=I\ts v\};
\een
cf. \cite[Sec.~3.3]{f:lc} and \cite[Sec.~8]{fr:qc}.
One can regard $\z_q(\wh\gl_n)$ as a subspace of $\Y_q(\gl_n)$.
Moreover, this subspace is closed under the multiplication in
the quantum affine algebra. Therefore, $\z_q(\wh\gl_n)$ can be identified
with a subalgebra of $\Y_q(\gl_n)$.
For $k=1,\dots,n$ introduce the series $\bar\ell_k(z)$
with coefficients in $\Y_q(\gl_n)$ by
\beql{barell}
\bar\ell_k(z)=\tr^{}_{1,\dots,k}\ts A^{(k)}\tss L^+_1(z)\dots L^+_k(zq^{-2k+2})
D_1\dots D_k.
\eeq

\bco\label{cor:inv}
All coefficients of the series
$\bar\ell_k(z)\vac$
with $k=1,\dots,n$
belong to the algebra of invariants $\z_q(\wh\gl_n)$.
Moreover, the coefficients of all series $\bar\ell_k(z)$ pairwise commute.
\eco

\bpf
By Theorem~\ref{thm:centr}, we have $L^-(u)\tss \ell_k(z)=\ell_k(z)\tss L^-(u)$.
Apply both sides to the vector $\vac\in V_q(\gl_n)$
and observe that $\ell_k(z)\vac=\bar\ell_k(z)\vac$.
This proves the first part of the corollary. The second part
follows by the application of both sides of the identity
$\ell_k(z)\ell_m(w)=\ell_m(w)\ell_k(z)$ to the vector $\vac$.
For the left hand side we get
\ben
\ell_k(z)\ell_m(w)\tss\vac=\ell_k(z)\bar\ell_m(w)\tss\vac
=\bar\ell_m(w)\tss \ell_k(z)\tss\vac=\bar\ell_m(w)\tss \bar\ell_k(z)\tss\vac.
\een
The same calculation for the right hand side gives
$\bar\ell_k(z)\bar\ell_m(w)=\bar\ell_m(w)\bar\ell_k(z)$.
\epf

The second part of the corollary is well known;
the series $\bar\ell_k(z)$ essentially coincides with the {\em transfer matrix}
associated with the $k$-th fundamental representation of $\U_q(\wh\gl_n)$;
see e.g. \cite{fr:qc}. The Harish-Chandra image of $\bar\ell_k(z)$
coincides with the $q$-{\em character} of this representation; see also
Theorem~\ref{thm:hch} below which recovers the calculation of the image
in a more general context.

\bre\label{rem:qva}
The form of the series $\ell_k(z)$ and $\bar\ell_k(z)$
indicates a possible interpretation of their properties from the viewpoint of the
quantum vertex algebra theory of \cite{ek:ql}.
\qed
\ere

\section{Quantum minor formulas for $\ell_k(z)$ and $\bar\ell_k(z)$}
\label{sec:qmf}
\setcounter{equation}{0}

By calculating the trace in \eqref{lkz}, we can get two
quantum minor-type expressions for $\ell_k(z)$
in terms of the entries of the matrices $L^{+}(z)=\big[l^+_{ij}(z)\big]$
and $\wt L(z):=L^-(z)^{-1}D=\big[\tss\wt l_{ij}(z)\big]$.
We will denote by
$l(\si)$ the length of a reduced decomposition of a permutation $\si\in\Sym_k$.
The length $l(\si)$ coincides with the number of inversions
in the sequence $\big(\si(1),\dots,\si(k)\big)$.

\bpr\label{prop:explk}
For $k=1,\dots,n$ we have
\begin{multline}
\ell_k(z)=\sum_{j_1,\dots,j_k}\ts \sum_{i_1<\dots<i_k}\ts\sum_{\si\in\Sym_k}
(-q)^{-l(\si)}\ts l^+_{i_{\si(1)}j_1}(z)\dots l^+_{i_{\si(k)}j_k}(z\tss q^{-2k+2})\\
{}\times\ts\wt l_{j_ki_k}(z\tss q^{-n-2k+2})\dots \wt l_{j_1i_1}(z\tss q^{-n})
\label{lkzlow}
\end{multline}
and
\begin{multline}
\ell_k(z)=\sum_{j_1,\dots,j_k}\ts \sum_{i_1<\dots<i_k}\ts\sum_{\si\in\Sym_k}
(-q)^{l(\si)}\ts l^+_{i_{\si(k)}j_k}(z)\dots l^+_{i_{\si(1)}j_1}(z\tss q^{-2k+2})\\
{}\times\ts\wt l_{j_1i_1}(z\tss q^{-n-2k+2})\dots \wt l_{j_ki_k}(z\tss q^{-n}).
\label{lkzhigh}
\end{multline}
\epr

\bpf
Using \eqref{lkz}
we interpret
\ben
A^{(k)}
L^+_1(v_1)\dots L^+_k(v_k)
\wt L_k(v_kq^{-n})\dots \wt L_1(v_1q^{-n})
\een
with $v_a=zq^{-2a+2}$
as an operator in the vector space $(\CC^n)^{\ot k}$.
Since
it is divisible on the right by $A^{(k)}$, the trace of the operator can be found
as $k!$ times the sum of the diagonal matrix elements corresponding to basis vectors
of the form $e_{i_1}\ot\dots\ot e_{i_k}$ with $i_1<\dots<i_k$. Now
\eqref{lkzlow} follows with the use of \eqref{antisym}
and the action of the $q$-permutations
on the basis vectors of this form: for any $\si\in\Sym_k$
\ben
P^{\tss q}_{\si}\big(e_{i_1}\ot\dots\ot e_{i_k}\big)
=q^{l(\si)}
\ts e_{i_{\si^{-1}(1)}}\ot\dots\ot e_{i_{\si^{-1}(k)}}
\een
and hence
\ben
P^{\tss q}_{\si}\big(e_{i_{\si(1)}}\ot\dots\ot e_{i_{\si(k)}}\big)=q^{-l(\si)}
\ts e_{i_1}\ot\dots\ot e_{i_k}.
\een
The proof of \eqref{lkzhigh} is quite similar;
use the basis vectors $e_{i_k}\ot\dots\ot e_{i_1}$
with the same condition $i_1<\dots<i_k$ on the indices.
\epf

\bre\label{rem:twomore}
Two more formulas for $\ell_k(z)$ analogous to \eqref{lkzlow} and \eqref{lkzhigh}
can be obtained by using
\eqref{lkzinv} instead of \eqref{lkz}.
\qed
\ere

The series
$\ell_n(z)$ can be factorized into a product of two quantum determinants.
To derive the factorization formula,
recall a construction of the quantum minors
of the matrices $L^{\pm}(z)$. Lemma~\ref{lem:deca} implies the relations
\beql{anitt}
A^{(k)} L^{\pm}_1(z)\dots L^{\pm}_k(q^{-2k+2}z)=L^{\pm}_k(q^{-2k+2}z)\dots L^{\pm}_1(z)  A^{(k)}.
\eeq
The {\it quantum minors\/} ${L^{\pm}(z)\tss}^{a_1\cdots\ts a_k}_{b_1\cdots\ts b_k}$
are the coefficients in the expansion of the either side of \eqref{anitt}
along the basis of matrix units:
\beql{matelmi}
\sum_{a_i,b_i}e_{a_1b_1}\ot\dots\ot e_{a_kb_k}
\ot {L^{\pm}(z)\tss}^{a_1\dots\ts a_k}_{b_1\dots\ts b_k}.
\eeq
The following formulas are immediate from the definition.
If $a_1<\dots<a_k$ then
\beql{qminorgenl}
{L^{\pm}(z)\tss}^{a_1\cdots\ts a_k}_{b_1\cdots\ts b_k}=
\sum_{\si\in \Sym_k} (-q)^{-l(\si)} \ts l^{\pm}_{a_{\si(1)}b_1}(z)\cdots
l^{\pm}_{a_{\si(k)}b_k}(q^{-2k+2}z)
\eeq
and for any $\tau\in\Sym_k$ we have
\beql{qmsym}
{L^{\pm}(z)\tss}^{a_{\tau(1)}\dots\ts a_{\tau(k)}}_{b_1\dots\ts b_k}=
(-q)^{l(\tau)}{L^{\pm}(z)\tss}^{a_1\dots\ts a_k}_{b_1\dots\ts b_k}.
\eeq
If $b_1<\dots<b_k$ (and the $a_i$ are arbitrary) then
\beql{qminorgen2}
{L^{\pm}(z)\tss}^{a_1\dots\ts a_k}_{b_1\dots\ts b_k}=
\sum_{\si\in \Sym_k} (-q)^{l(\si)} \ts l^{\pm}_{a_kb_{\si(k)}}(q^{-2k+2}z)\dots
l^{\pm}_{a_1b_{\si(1)}}(z),
\eeq
and for any $\tau\in\Sym_k$ we have
\beql{qmsym2}
{L^{\pm}(z)\tss}^{a_1\dots\ts a_k}_{b_{\tau(1)}\dots\ts b_{\tau(k)}}=
(-q)^{-l(\tau)}{L^{\pm}(z)\tss}^{a_1\dots\ts a_k}_{b_1\dots\ts b_k}.
\eeq
Moreover, the quantum minor is zero if two top or two bottom indices
are equal.

The following lemma is well-known. We give a
proof for completeness.\footnote{The critical level assumption $k=-2$
is omitted in the corresponding statement in \cite[Lemma~2]{fr:qa} in the case $n=2$.
It was only used there under this assumption.}

\ble\label{lem:qd}
The coefficients of the quantum determinants
\beql{qdet}
\qdet L^{\pm}(z)={L^{\pm}(z)\tss}^{1\cdots\ts n}_{1\cdots\ts n}
\eeq
belong to the center of the quantum affine algebra $\U_q(\wh\gl_n)_{\text{\rm cri}}$
at the critical level.
\ele

\bpf
Introduce the product
\ben
R(v_0,v_1,\dots,v_n)=\prod_{0\leqslant a<b\leqslant n}R_{ab}(v_a/v_b),
\een
where the $v_a$ are variables and the product is taken in the
lexicographical order on the pairs $(a,b)$. The defining relations
\eqref{RLL} and \eqref{RLLpm} imply
\ben
R(uq^n,v_1,\dots,v_n)L_0^{+}(u)L_1^{-}(v_1)\dots L_n^{-}(v_n)
=L_n^{-}(v_n)\dots L_1^{-}(v_1)L_0^{+}(u)R(uq^{-n},v_1,\dots,v_n).
\een
Use \eqref{ruv}--\eqref{rf} to write this relation in terms of the $R$-matrix $R(u,v)$.
By cancelling common factors we get
\ben
\wt R(uq^n,v_1,\dots,v_n)L_0^{+}(u)L_1^{-}(v_1)\dots L_n^{-}(v_n)
=L_n^{-}(v_n)\dots L_1^{-}(v_1)L_0^{+}(u)\wt R(uq^{-n},v_1,\dots,v_n),
\een
where
\ben
\wt R(v_0,v_1,\dots,v_n)=\prod_{a=1}^n \frac{f(v_0/v_a)}{\qin v_0-q v_a}\ts
\prod_{0\leqslant a<b\leqslant n}R_{ab}(v_a,v_b).
\een
Now specialize the variables by setting $v_a=zq^{-2a+2}$ for $a=1,\dots,n$ and replace
the product of $R$-matrices $R_{ab}(v_a,v_b)$
over the set of pairs $1\leqslant a<b\leqslant n$ using Lemma~\ref{lem:deca}.
Since
\beql{ald}
A^{(n)} L^{\pm}_1(v_1)\dots L^{\pm}_n(v_n)=L^{\pm}_n(v_n)\dots L^{\pm}_1(v_1)\tss A^{(n)}=
A^{(n)}\qdet L^{\pm}(z),
\eeq
we get
\begin{multline}
\prod_{a=1}^n \frac{f(uq^n/v_a)}{u\tss q^{n-1}-q v_a}
\ts\prod_{a=1,\dots,n}^{\longrightarrow}
R_{0a}(uq^n,v_a)
A^{(n)}L_0^{+}(u)\ts\qdet L^{-}(z)\\
{}=\prod_{a=1}^n \frac{f(uq^{-n}/v_a)}{u\tss q^{-n-1}-q v_a}
\ts\qdet L^{-}(z)\tss L_0^{+}(u)\tss A^{(n)}
\ts\prod_{a=1,\dots,n}^{\longleftarrow}  R_{0a}(uq^{-n},v_a).
\non
\end{multline}
Observe that
\ben
\prod_{a=1,\dots,n}^{\longrightarrow}
R_{0a}(v_0,v_a)A^{(n)}=
A^{(n)}
\ts\prod_{a=1,\dots,n}^{\longleftarrow}  R_{0a}(v_0,v_a)
=A^{(n)}\tss(\qin v_0-qv_1)\prod_{a=2}^n (v_0-v_a).
\een
Indeed, by the first equality, it suffices to verify the second equality
on the basis vectors of the form
$e_i\ot e_i\ot e_1\ot\dots\ot e_{i-1}\ot e_{i-1}\ot\dots\ot e_n$
for $i=1,\dots,n$ which is straightforward. Thus, we can conclude that
$L_0^{+}(u)$ commutes with $\qdet L^{-}(z)$ due to the identity
\ben
\prod_{a=1}^n \frac{f(uq^n/v_a)}{f(uq^{-n}/v_a)}
=\prod_{a=2}^n \frac{(u\tss q^{n-1}-q v_a)(u\tss q^{-n}-v_a)}
{(u\tss q^{-n-1}-q v_a)(u\tss q^{n}-v_a)}
\een
which follows from \eqref{fx}. The relation
$L_0^{-}(u)\ts\qdet L^{+}(z)=\qdet L^{+}(z)L_0^{-}(u)$
is verified by a similar argument with the use of the unitarity property
of the $R$-matrix \eqref{rx}:
\ben
\overline R(x^{-1})=\overline R_{21}(x)^{-1}.
\een
The proof of the remaining two relations
$L_0^{\pm}(u)\ts\qdet L^{\pm}(z)=\qdet L^{\pm}(z)L_0^{\pm}(u)$
is simpler as it relies only on the defining relations \eqref{RLL}.
\epf

\bre\label{rem:levelzero}
Both quantum determinants \eqref{qdet} are also known to be central in
the quantum affine algebra $\U_q(\wh\gl_n)$ at the zero level $c=0$; the algebra is
defined as in Sec.~\ref{sec:qaa} with the factor $f(x)$ in \eqref{rx} omitted.
\qed
\ere

\bco\label{cor:qdets}
We have
\ben
\ell_n(z)=\qdet L^+(z)\ts \big(\qdet L^-(zq^{-n})\big)^{-1}.
\een
\eco

\bpf
Relation \eqref{ald} implies
\ben
A^{(n)}L^-_n(v_n)^{-1}\dots L^-_1(v_1)^{-1} =\big(\qdet L^-(z)\big)^{-1}\tss A^{(n)}.
\een
Replacing here $z$ by $zq^{-n}$ and using $\qdet D=1$ we get the desired formula
from \eqref{lkz}.
\epf

We will now give a formulation of Corollary~\ref{cor:inv}
by combining the series $\bar\ell_k(z)$ defined in \eqref{barell} into a single determinant
by some analogy with \cite{cm:ho} and \cite{ct:qs}.
Introduce the extension $\Y^{\text{\rm ext}}_q(\gl_n)$ of the algebra
$\Y_q(\gl_n)$ by adjoining pairwise commuting elements $\pi_1,\dots,\pi_n$
subject to the additional relations
\beql{pilkm}
\pi_i\ts l^+_{km}(u)=\begin{cases} l^+_{km}(u)\ts \pi_i\qquad&\text{if}\quad i\leqslant k,\\
q^2\tss l^+_{km}(u)\ts \pi_i\qquad&\text{if}\quad i> k.
\end{cases}
\eeq
Combine the elements $\pi_i$ into the diagonal matrix $\Pi=\diag\ts[\pi_1,\dots,\pi_n]$.
We have a vector space isomorphism $\Y^{\text{\rm ext}}_q(\gl_n)/J\cong \Y_q(\gl_n)$,
where $J$ is the left ideal of $\Y^{\text{\rm ext}}_q(\gl_n)$ generated
by the elements $\pi_i-1$ with $i=1,\dots,n$. We will identify $\Y_q(\gl_n)$
with the quotient via this isomorphism.

We will point out a connection with
$q$-analogues of Manin matrices (also known as {\em right quantum matrices}); see, e.g.,
\cite{cfrs:ap} for a detailed account of their properties.
An $n\times n$ matrix $M$ with entries in an associative algebra $\Ac$
is called $q$-{\em Manin} if it satisfies the relation
\ben
A^{(2)}M_1M_2=A^{(2)}M_1M_2A^{(2)}
\een
in the algebra $\End \CC^n\ot\End \CC^n\ot \Ac$ with the meaning
of the subscripts as in \eqref{la}.

Introduce the operator $\de$ which interacts with power series in $z$ by the rule
$\de\tss g(z)=g(zq^{-2})\tss\de$. Adjoining this element to the algebra $\Y_q(\gl_n)$
we find that both $L^+(z)\tss\de$ and $L^+(z)\tss D\tss\de$ are $q$-Manin matrices.
Define the $q$-{\em determinant} of a square matrix $M$ by
\ben
{\det}_q\tss M=\sum_{\si\in \Sym_n} (-q)^{-l(\si)} \ts M_{\si(1)1}\cdots
M_{\si(n)n}.
\een
In the following proposition we regard the $q$-determinant of the matrix $\Pi+L^+(z)\tss D\tss\de$
as a polynomial in $\de$ with coefficients in $\Y^{\text{\rm ext}}_q(\gl_n)$.

\bpr\label{prop:detq}
We have the relation modulo the left ideal $J$:
\ben
{\det}_q\big(\Pi+L^+(z)\tss D\tss\de\big)=1+\sum_{k=1}^n \bar\ell_k(z)\ts \de^{\tss k}.
\een
\epr

\bpf
Calculating the $q$-determinant
we will write it as the sum of monomials of the form
\begin{multline}
(-q)^{-l(\si)} M_{\si(1)\tss 1}\cdots M_{\si(i_1-1)\ts i_1-1}\tss \pi_{i_1}
M_{\si(i_1+1)\ts i_1+1}\\[0.4em]
{}\times\cdots M_{\si(i_2-1)\ts i_2-1}\tss \pi_{i_2}
M_{\si(i_2+1)\ts i_2+1}\cdots \pi_{i_{n-k}}
M_{\si(i_{n-k}+1)\ts i_{n-k}+1}\cdots M_{\si(n)n},
\non
\end{multline}
where $M=L^+(z)\tss D\tss\de$ and $\si(i_a)=i_a$ for $a=1,\dots,n-k$.
Now use relations \eqref{pilkm} to move all the elements $\pi_{i_a}$
so they will appear to the right from all factors $M_{km}$
in the monomial. As a result, by moving each element $\pi_{i_a}$ we get
the factor $q^{2r_a}$, where $r_a$ is the number
of elements of the set $\{\si(i_a+1),\dots,\si(n)\}$ which are less than $i_a$.
That is, $r_a$ is the number of inversions
formed by the index $i_a$ with the indices of the set $\{\si(i_a+1),\dots,\si(n)\}$.
Therefore, after moving all the elements $\pi_{i_a}$,
the monomial will get
the factor $(-q)^{-l(\wh\si)}$, where $\wh\si$ is the sequence
of elements obtained from $(\si(1),\dots,\si(n))$ by removing $i_1,\dots,i_{n-k}$,
and $l(\wh\si)$ is the number of inversions in that sequence.
This demonstrates that the coefficient of the product $\pi_{i_1} \dots \pi_{i_{n-k}}$
coincides with the $q$-determinant of the principal submatrix of $M$ obtained
by deleting rows and columns enumerated by $i_1,\dots,i_{n-k}$.
This $q$-determinant equals
$
\big[{L^{+}(z) D\tss}\big]^{a_1\cdots\ts a_k}_{a_1\cdots\ts a_k}\ts \de^k
$
where $\{a_1,\dots,a_k\}=\{1,\dots,n\}\setminus\{i_1,\dots,i_{n-k}\}$.
Since
\ben
\bar\ell_k(z)=\sum_{1\leqslant a_1<\dots<a_k\leqslant n}
\big[{L^{+}(z) D\tss}\big]^{a_1\cdots\ts a_k}_{a_1\cdots\ts a_k},
\een
the required relation follows by taking the quotient over the left ideal $J$.
\epf

The next lemma will be used in
Sec.~\ref{sec:hch} below.

\ble\label{lem:com}
The entries of the inverse matrix $L^-(z)^{-1}$ are found by the formula
\ben
\big[L^-(z)^{-1}\big]_{ij}=(-q)^{j-i}\ts \big(\qdet L^-(zq^{2n-2})\big)^{-1}\ts
{L^{-}(zq^{2n-2})\tss}^{1\dots\ts\wh j\dots n}_{1\dots\ts\wh i\dots n},
\een
where the hats indicate indices to be skipped.
\ele

\bpf
By \eqref{ald} we have
\ben
A^{(n)} L^{-}_1(v_1)\dots L^{-}_{n-1}(v_{n-1})=
A^{(n)}\qdet L^{-}(z)L^{-}_{n}(v_n)^{-1},
\een
where $v_a=zq^{-2a+2}$, as before.
The desired formula follows by the application of both sides to the basis
vector $e_1\ot\dots\ot\wh{e_i}\ot\dots\ot e_n\ot e_j$
and the replacement $z\mapsto z\tss q^{2n-2}$.
\epf

\section{Poincar\'e--Birkhoff--Witt theorem}
\label{sec:pbw}
\setcounter{equation}{0}

To define analogues of the Harish-Chandra homomorphism,
we will need a version of the Poincar\'e--Birkhoff--Witt theorem
for the quantum affine algebra $\U_q(\wh\gl_n)$.
Introduce a total ordering $\prec$ on the set of generators as follows.
First, each generator $l^+_{ij}[r]$ precedes each generator $l^-_{km}[s]$.
Furthermore, $l^+_{ij}[r]\prec l^+_{km}[s]$ if and only if the triple $(j-i,i,r)$
precedes $(m-k,k,s)$ in the lexicographical order. Finally, we set
$l^-_{ij}[r]\prec l^-_{km}[s]$ if and only if the triple $(i-j,i,r)$
precedes $(k-m,k,s)$ in the lexicographical order.
Note that by the defining relations \eqref{RLL},
\beql{commll}
\big[l^{\pm}_{ij}[r],l^{\pm}_{ij}[s]\big]=0
\eeq
for all $r$ and $s$. Hence, the ordering $\prec$ induces a well-defined total ordering
on the series \eqref{serlpm} such that $l^+_{ij}(u)\prec l^-_{km}(u)$ and
\ben
\bal
&l^+_{n\tss 1}(u)\prec l^+_{n-1\ts 1}(u)\prec l^+_{n\ts 2}(u)\prec\dots\prec
l^+_{1\tss 1}(u)\prec\dots\prec
l^+_{n\tss n}(u)\prec l^+_{1\tss 2}(u)\prec\dots\prec l^+_{1\tss n}(u),\\[0.3em]
&l^-_{1\tss n}(u)\prec l^-_{1\ts n-1}(u)\prec l^-_{2\ts n}(u)\prec\dots\prec
l^-_{1\tss 1}(u)\prec\dots\prec
l^-_{n\tss n}(u)\prec l^-_{2\tss 1}(u)\prec\dots\prec l^-_{n\tss 1}(u).
\eal
\een
Consider the ordered monomials in the generators $l^{\pm}_{ij}[r]$
multiplied by integer powers of the central element $\ga=q^c$
(the zero elements
$l^+_{ij}[0]$ for $i>j$
and $l^-_{ij}[0]$ for $i<j$ are excluded).
Relations \eqref{RLL} and \eqref{RLLpm} imply that
\beql{lpze}
l^+_{ii}[0]\tss l^{\pm}_{km}(u)=q^{-\de_{ik}+\de_{im}}\tss l^{\pm}_{km}(u)\tss l^+_{ii}[0]
\Fand
l^-_{ii}[0]\tss l^{\pm}_{km}(u)=q^{\de_{ik}-\de_{im}}\tss l^{\pm}_{km}(u)\tss l^-_{ii}[0].
\eeq
Hence,
using \eqref{lplm}
we may suppose that for each $i=1,\dots,n$
each monomial only contains either a nonnegative power of $l^+_{ii}[0]$
or a positive power of $l^-_{ii}[0]$.
Under these assumptions we have the following version of the
Poincar\'e--Birkhoff--Witt theorem.

\bpr\label{prop:pbw}
The ordered monomials in the generators form a basis of the
quantum affine algebra $\U_q(\wh\gl_n)$.
\epr

\bpf
First, we prove the claim for the quantum affine algebra $\overline \U_q(\wh\gl_n)$
which is defined in the same way as $\U_q(\wh\gl_n)$; the only difference
is the use of the $R$-matrix \eqref{rx} instead of \eqref{rf}.
Thus, we only replace \eqref{RLLpm} with the relation
\beql{RLLpmbar}
\overline R(u\ga^{-1}/v)\overline L_1^{\ts +}(u)\overline L_2^{\ts -}(v)
=\overline L_2^{\ts -}(v)\overline L_1^{\ts +}(u)\overline R(u\ga/v)
\eeq
and leave all other defining relations unchanged. Here we use the
bar symbol over the respective objects associated with
the algebra $\overline \U_q(\wh\gl_n)$.
We begin by showing that the ordered monomials in the generators of
$\overline \U_q(\wh\gl_n)$ span the algebra. Given a monomial
in the generators we will use the induction on its length to show that
it equals a linear combination of ordered monomials. Writing
\eqref{RLLpmbar} in terms of the entries of the matrices $\overline L^{\ts \pm}(u)$
we get
\begin{align}
(\ga\ts q^{-\delta_{jm}} u - q^{\delta_{jm}}v)\ts
\bar l_{km}^-(v)\ts& \bar l_{ij}^+(u)
+ (q^{-1} - q)\ts (u\ts\ga\ts \delta_{m>j} +v\ts \delta_{m<j})\ts
\bar l_{kj}^-(v)\ts \bar l_{im}^+(u)
\non\\[0.4em]
{}=\Big[(\ga^{-1}\ts q^{-\delta_{ik}} u - q^{\delta_{ik}}v)
\ts&\bar l_{ij}^+(u)\ts \bar l_{km}^-(v)
\non\\
{}+ (q^{-1} - q)\ts &(u\ts\ga^{-1}\ts \delta_{i>k} +v\ts \delta_{i<k})
\ts\bar l_{kj}^+(u)\ts \bar l_{im}^-(v)\Big]\ts
\frac{\ga\ts q^{-2} u - v}{\ga^{-1}\ts q^{-2} u - v},
\label{lplmga}
\end{align}
where $\de_{i<j}$ or $\de_{i>j}$ equals $1$ if
the subscript inequality
is satisfied and $0$ otherwise. If $m=j$ then the relation allows us to write
$\bar l_{kj}^-[s]\ts \bar l_{ij}^+[r]$ as a linear combination of ordered products
of generators. If $m\ne j$ then we swap $m$ and $j$ in the relation to get a system of two equations
for $\bar l_{km}^-(v)\ts \bar l_{ij}^+(u)$ and $\bar l_{kj}^-(v)\ts \bar l_{im}^+(u)$.
By solving the system we will be able to write $\bar l_{km}^-[s]\ts \bar l_{ij}^+[r]$
as a linear combination of ordered products of generators.

Similarly, \eqref{RLL} gives the relation
\begin{align}
&(q^{-\delta_{ik}} u - q^{\delta_{ik}}v)
\ts\bar l_{ij}^{\ts\pm}(u)\ts \bar l_{km}^{\ts\pm}(v)
+ (q^{-1} - q)\ts (u\ts \delta_{i>k} +v\ts \delta_{i<k})
\ts\bar l_{kj}^{\ts\pm}(u)\ts \bar l_{im}^{\ts\pm}(v)
\non\\[0.3em]
&{}=(q^{-\delta_{jm}} u - q^{\delta_{jm}}v)\ts
\bar l_{km}^{\ts\pm}(v)\ts \bar l_{ij}^{\ts\pm}(u)
+ (q^{-1} - q)\ts (u\ts \delta_{m>j} +v\ts \delta_{m<j})\ts
\bar l_{kj}^{\ts\pm}(v)\ts \bar l_{im}^{\ts\pm}(u)
\label{ll}
\end{align}
which implies that $\bar l_{ij}^{\ts\pm}[r]\tss\bar l_{km}^{\ts\pm}[s]$
is a linear combination of ordered products of generators; see \cite[Corollary~2.13]{gm:rt}
for a detailed argument.\footnote{The generator matrices $T(u)$ and $\overline T(u)$
of \cite{gm:rt} correspond to $\overline L^{\tss-}(u)$ and $\overline L^{\tss +}(u)$, respectively.}

As a next step, we will show that the ordered monomials are linearly
independent in $\overline \U_q(\wh\gl_n)$. Consider the quantum affine algebra
$\U_q(\wh\gl_n)_0$ at the level zero which is the quotient of $\overline \U_q(\wh\gl_n)$
by the relation $\ga=1$.
We have the natural epimorphism
\ben
\psi:\overline \U_q(\wh\gl_n)\to \U_q(\wh\gl_n)_0, \qquad
\overline L^{\ts \pm}(u)\mapsto L^{\ts \pm}(u),\qquad \ga\mapsto 1.
\een
Suppose that a Laurent polynomial in $\ga$, whose coefficients are
nontrivial linear combinations of
ordered monomials in the generators $\bar l_{ij}^{\ts\pm}[r]$, is zero in
$\overline \U_q(\wh\gl_n)$. Multiplying by a power of $\ga$ if necessary, we get
a polynomial in $\ga$ equal to zero. Choose such a polynomial
$P=x_k\tss\ga^k+\dots+x_0$
of the minimal possible degree $k\geqslant 0$.
Since $P=0$ in $\overline \U_q(\wh\gl_n)$ we have
\ben
0=\psi(P)=\psi(x_k)+\dots+\psi(x_0)=\psi(x_k+\dots+x_0).
\een
The sum $x=x_k+\dots+x_0$ is a linear combination of
ordered monomials in the generators $\bar l_{ij}^{\ts\pm}[r]$.
By the definition of $\psi$, the image $\psi(x)$ is
the corresponding linear combination of
ordered monomials in the generators $l_{ij}^{\ts\pm}[r]$ of $\U_q(\wh\gl_n)_0$.
On the other hand, by the arguments of \cite[Sec.~2.3]{gm:rt} applied to
this particular ordering, the corresponding version of the
Poincar\'e--Birkhoff--Witt theorem holds for $\U_q(\wh\gl_n)_0$, so that the
ordered monomials are linearly independent.
Hence, $\psi(x)=0$ implies $x=0$.
If $k=0$ then this is a contradiction. If $k\geqslant 1$ we can write
\ben
P=x_k\tss\ga^k+\dots+x_0=(\ga-1)(y_{k-1}\tss\ga^{k-1}+\dots+y_0),
\een
where the $y_i$ are again linear combinations of
ordered monomials in the generators $\bar l_{ij}^{\ts\pm}[r]$.
By the results of \cite{df:it}, the algebra $\overline \U_q(\wh\gl_n)$
can be defined by the Drinfeld generators. Due to the well-known
relationship between the quantum affine algebras associated with $\sll_n$ and $\gl_n$
(see, e.g., \cite[Sec.~2.6]{fm:ha}), the
Poincar\'e--Birkhoff--Witt theorem for the algebra $\U_q(\wh\sll_n)$ in its Drinfeld
presentation \cite{b:cb, b:bg} implies that the relation
$P=0$ is possible only if
$y_{k-1}\tss\ga^{k-1}+\dots+y_0=0$. This contradicts the minimality of
the degree $k$ thus completing the proof for $\overline \U_q(\wh\gl_n)$.

Finally, we extend the argument to the quantum affine algebra $\U_q(\wh\gl_n)$
defined with the $R$-matrix \eqref{rf} so that \eqref{RLLpm} with $\ga=q^c$ should be used
instead of \eqref{RLLpmbar}. This affects only relation \eqref{lplmga} (for the
generators $l_{ij}^{\ts\pm}[r]$ instead of $\bar l_{ij}^{\ts\pm}[r]$) which will now
get an extra factor $f(u\ga^{-1}/v)/f(u\ga/v)$ on the right hand side.
However, this does not bring any change into the first part of the argument showing that
the ordered monomials in the generators span the algebra $\U_q(\wh\gl_n)$.

To prove the linear independence of the ordered monomials, we follow \cite[Sec.~V]{df:it}
and introduce a homomorphism
\beql{homheis}
\phi:\U_q(\wh\gl_n)\to\Hc_q(n)\ot_{\CC[\ga,\ga^{-1}]} \overline \U_q(\wh\gl_n),
\eeq
where $\Hc_q(n)$ is the Heisenberg algebra with generators $\ga$ and $h[r]$, $r\in\ZZ$, $r\ne 0$.
The defining relations of $\Hc_q(n)$ have the form
\ben
\big[h[r],h[s]\big]=\de_{r, -s}\ts \al[r],\qquad r\geqslant 1,
\een
and $\ga$ is central and invertible; all other pairs of the generators commute.
The elements $\al[r]$ are defined by the expansion
\ben
\exp\ts\sum_{r=1}^{\infty}\al[r]\ts x^r=\frac{f(x\ga)}{f(x\ga^{-1})}.
\een
So we have the identity
\begin{multline}
f(u\ga^{-1}/v)\ts \exp\ts\Big(\sum_{r=1}^{\infty}h[r]\ts u^r\Big)
\ts \exp\ts\Big(\sum_{s=1}^{\infty}h[-s]\ts v^{-s}\Big)\\
{}=f(u\ga/v)\ts \exp\ts\Big(\sum_{s=1}^{\infty}h[-s]\ts v^{-s}\Big)
\ts \exp\ts\Big(\sum_{r=1}^{\infty}h[r]\ts u^r\Big).
\non
\end{multline}
Clearly, the monomials of the form
$h[r_1]\dots h[r_k]$ with $k\geqslant 0$ and $r_1\geqslant\dots\geqslant r_k$
(with $r_i\ne 0$)
form a basis of the $\CC[\ga,\ga^{-1}]$-module $\Hc_q(n)$.
The homomorphism \eqref{homheis} is now defined by $\phi:\ga\mapsto\ga$ and
\beql{heisext}
\phi:L^+(u)\mapsto \exp\ts\Big(\sum_{r=1}^{\infty}h[r]\ts u^r\Big)\ts
\overline L^{\ts +}(u),\quad
L^-(u)\mapsto \exp\ts\Big(\sum_{r=1}^{\infty}h[-r]\ts u^{-r}\Big)\ts
\overline L^{\ts -}(u).
\eeq
Suppose there is a linear combination of the ordered monomials in the
generators of $\U_q(\wh\gl_n)$ equal to zero. Consider its image under
the homomorphism $\phi$. Using the basis $\{h[r_1]\dots h[r_k]\}$
of $\Hc_q(n)$ and the Poincar\'e--Birkhoff--Witt basis for the algebra $\overline \U_q(\wh\gl_n)$,
we conclude that all coefficients of the linear combination
must be zero.
\epf

We will also need a version of the Poincar\'e--Birkhoff--Witt theorem
for a different ordering of the generators. As we pointed out above,
a total ordering can be defined on the generating series \eqref{serlpm}
due to \eqref{commll}. We set $l^+_{ij}(u)\prec l^-_{km}(u)$ as before,
but the remaining conditions are swapped between $l^+_{ij}(u)$ and $l^-_{ij}(u)$:
\ben
\bal
&l^+_{1\tss n}(u)\prec l^+_{1\ts n-1}(u)\prec l^+_{2\ts n}(u)\prec\dots\prec
l^+_{1\tss 1}(u)\prec\dots\prec
l^+_{n\tss n}(u)\prec l^+_{2\tss 1}(u)\prec\dots\prec l^+_{n\tss 1}(u),\\[0.3em]
&l^-_{n\tss 1}(u)\prec l^-_{n-1\ts 1}(u)\prec l^-_{n\ts 2}(u)\prec\dots\prec
l^-_{1\tss 1}(u)\prec\dots\prec
l^-_{n\tss n}(u)\prec l^-_{1\tss 2}(u)\prec\dots\prec l^-_{1\tss n}(u).
\eal
\een
Under the same assumptions on the monomials as for Proposition~\ref{prop:pbw},
the following holds.

\bpr\label{prop:pbwopp}
The ordered monomials in the generators form a basis of the
quantum affine algebra $\U_q(\wh\gl_n)$.
\epr

\bpf
The argument is the same as for Proposition~\ref{prop:pbw} with some obvious minor
changes taking into the account the ordering conditions.
\epf

\section{Harish-Chandra homomorphisms}
\label{sec:hch}
\setcounter{equation}{0}

Consider the quantum affine algebra $\U_q(\wh\gl_n)_{\text{\rm cri}}$
at the critical level, $\ga=q^{-n}$.
By Proposition~\ref{prop:pbw}, any element $x\in\U_q(\wh\gl_n)_{\text{\rm cri}}$
can be written as a unique linear combination of
ordered monomials
in the generators $l^{\pm}_{ij}[r]$. Denote by $\U^0$
the subspace of the algebra spanned by
those monomials which do not
contain any generators $l^{\pm}_{ij}[r]$ with $i\ne j$.
Let $x_0$ denote the component
of the linear combination representing the element $x$,
which belongs to $\U^0$.
The mapping $\theta:x\mapsto x_0$ defines the projection
$\theta:\U_q(\wh\gl_n)_{\text{\rm cri}}\to \U^0$. Extending it
by continuity we get the projection $\theta:\wt\U_q(\wh\gl_n)_{\text{\rm cri}}\to \wt\U^0$
to the corresponding completed vector space $\wt\U^0$.

Introduce the algebra $\Pi_q(n)$ as
the quotient of the algebra of polynomials in independent variables
$l^+_i[-r]$, $l^-_i[r]$ with $i=1,\dots,n$ and $r=0,1,\dots$
by the relations $l^+_i[0]\tss l^-_i[0]=1$ for all $i$.
The mapping $\eta:\U^0\to \Pi_q(n)$ which takes each ordered monomial
in the generators $l^{\pm}_{ii}[\mp r]$ to the corresponding monomial
in the variables $l^{\pm}_{i}[\mp r]$ by the rule
$l^{\pm}_{ii}[\mp r]\mapsto l^{\pm}_{i}[\mp r]$ extends to an isomorphism
of vector spaces.
Define
the completion $\wt\Pi_q(n)$ of the algebra $\Pi_q(n)$ as
the inverse limit
\ben
\wt\Pi_q(n)=\lim_{\longleftarrow} \Pi_q(n)/I_p, \qquad p>0,
\een
where $I_p$ denotes the ideal of $\Pi_q(n)$ generated by all elements
$l^{-}_{i}[r]$ with $r\geqslant p$; cf. \eqref{compl}. The isomorphism
$\eta$ extends to an isomorphism of the respective completed
vector spaces $\eta:\wt\U^0\to \wt\Pi_q(n)$.
Thus we get a linear map
\beql{chihom}
\chi:\wt\U_q(\wh\gl_n)_{\text{\rm cri}}\to \wt\Pi_q(n)
\eeq
defined as the composition $\chi=\eta\circ\theta$.
The next proposition provides an analogue of the Harish-Chandra homomorphism
for the quantum affine algebra.

\bpr\label{prop:hchhom}
The restriction of the map \eqref{chihom} to the center $\Zr_q(\wh\gl_n)$ of the algebra
$\wt\U_q(\wh\gl_n)_{\text{\rm cri}}$ is a homomorphism of commutative algebras
\beql{hch}
\chi:\Zr_q(\wh\gl_n)\to \wt\Pi_q(n).
\eeq
\epr

\bpf
For $x,y\in \Zr_q(\wh\gl_n)$ set $x_0=\chi(x)$ and $y_0=\chi(y)$.
Write $y$ as a (possibly infinite) linear combination of ordered monomials
in the generators $l^{\pm}_{ij}[r]$. Suppose that
\ben
m=\prod_a l^{+}_{i_aj_a}[r_a]\ts\prod_b l^{-}_{i_bj_b}[r_b]
\een
is an ordered monomial which occurs in the linear combination.
Note its property
\beql{wt}
\sum_a(i_a-j_a)+\sum_b(i_b-j_b)=0
\eeq
implied by \eqref{lpze}. Suppose that
$m\in\ker\chi$. Since $x$ is in the center, we have
\ben
x\tss m=\prod_a l^{+}_{i_aj_a}[r_a]\ts x\ts \prod_b l^{-}_{i_bj_b}[r_b].
\een
To write $x\tss m$ as a linear combination of ordered monomials
we will only need to use the defining relations \eqref{RLL}
which are also given in \eqref{ll} where the series $\bar l_{ij}^{\ts\pm}(u)$
should be replaced with $l_{ij}^{\ts\pm}(u)$, respectively.
Since the relations \eqref{ll} are homogeneous with respect to the weight parameter
$i-j+k-m$, we derive that
$x\tss m\in\ker\chi$. Hence a nonzero contribution to the image $\chi(xy)$
can only come from $\chi(xy_0)$, that is, from expressions of the form
\ben
\prod_a l^{+}_{i_ai_a}[r_a]\ts x\ts \prod_b l^{-}_{i_bi_b}[r_b].
\een
If $p$ is an ordered
monomial which occurs in the linear combination representing $x$ and $\chi(p)=0$,
then applying property \eqref{wt} to the monomial $p$ we conclude that
\ben
\chi:\prod_a l^{+}_{i_ai_a}[r_a]\ts p\ts \prod_b l^{-}_{i_bi_b}[r_b]\mapsto 0.
\een
Finally, observe that by the defining relations \eqref{ll},
any two generators
$l^{-}_{ii}[r]$ and $l^{-}_{jj}[s]$ (resp., $l^{+}_{ii}[r]$ and $l^{+}_{jj}[s]$)
can be permuted modulo $\ker\chi$ within any monomial of the
form
\ben
\prod_a l^{+}_{i_ai_a}[r_a]\ts \prod_b l^{-}_{i_bi_b}[r_b].
\een
This proves that $\chi(xy)=x_0y_0$.
\epf

Now we are in a position to calculate the Harish-Chandra images of the
higher Sugawara operators provided by Theorem~\ref{thm:centr}.
Combine the generators of the algebra $\Pi_q(n)$ into the series
\ben
l^{+}_{i}(z)=\sum_{r=0}^{\infty}l^{+}_{i}[-r]\tss z^r,\qquad
l^{-}_{i}(z)=\sum_{r=0}^{\infty}l^{-}_{i}[r]\tss z^{-r}
\een
and for $i=1,\dots,n$ set
\ben
\la_i(z)=q^{n-2i+1}\ts \frac{l^{+}_{i}(z)\ts l^{-}_{1}(zq^{-n+2})\dots l^{-}_{i-1}(zq^{-n+2i-2})}
{l^{-}_{1}(zq^{-n})\dots l^{-}_{i}(zq^{-n+2i-2})}.
\een
This is a Laurent series in $z$ whose coefficients
are elements of the completed algebra $\wt\Pi_q(n)$.

\bth\label{thm:hch}
For each $k=1,\dots,n$ the image of the series $\ell_k(z)$ under the
Harish-Chandra homomorphism \eqref{hch} is found by
\ben
\chi:\ell_k(z)\mapsto \sum_{1\leqslant i_1<\dots<i_k\leqslant n}
\la_{i_1}(z)\la_{i_2}(zq^{-2})\dots \la_{i_k}(zq^{-2k+2}).
\een
\eth

\bpf
We will use formula \eqref{lkzlow} for $\ell_k(z)$.
Apply Lemma~\ref{lem:com} to express the series $\wt l_{ji}(z)=q^{n-2i+1}\ts\big[L^-(z)^{-1}\big]_{ji}$
in terms of quantum minors. By \eqref{qmsym2} we have
\ben
{L^{-}(z)\tss}^{1\dots\ts\wh i\dots n}_{1\dots\ts\wh j\dots n}
=(-q)^{l(\om)}\ts {L^{-}(z)\tss}^{1\dots\ts\wh i\dots n}_{n\dots\ts\wh j\dots 1},
\een
where $\om\in\Sym_{n-1}$ reverses the order of the lower indices.
Expanding this quantum minor
by \eqref{qminorgenl} we find that
a nonzero contribution to the image $\chi\big(\ell_k(z)\big)$ can only come from
the summands in \eqref{lkzlow} with $i_{\si(1)}\leqslant j_1\leqslant i_1$.
These conditions imply that $\si(1)=1$ and $j_1=i_1$. By the defining relations
in $\U_q(\wh\gl_n)$, the same observation gives $\si(2)=2$ and $j_2=i_2$, etc.,
so that a nonzero contribution comes only from the terms with $\si=1$ and $j_a=i_a$ for
all $a=1,\dots,k$. Applying Lemma~\ref{lem:com} and formulas \eqref{qminorgenl} and \eqref{qmsym2} again
we find that the contributions of the quantum minors are found by
\ben
\qdet L^-(z)\mapsto l^-_1(zq^{-2n+2})\dots l^-_n(z)
\een
and
\ben
{L^{-}(z)\tss}^{1\dots\ts\wh i\dots n}_{1\dots\ts\wh i\dots n}
\mapsto l^-_1(zq^{-2n+4})\dots l^-_{i-1}(zq^{-2n+2i})\ts l^-_{i+1}(zq^{-2n+2i+2})\dots l^-_n(z).
\een
This completes
the calculation of the Harish-Chandra image of $\ell_k(z)$.
\epf

Consider the restriction of
the map \eqref{chihom} to the subalgebra $\Y_q(\gl_n)$ of $\wt\U_q(\wh\gl_n)_{\text{\rm cri}}$.
As in Sec.~\ref{sec:mt}, we impose the conditions $l^+_{ii}[0]=1$ for all $i$
so that
\beql{chihomy}
\chi:\Y_q(\gl_n)\to \Pi^+_q(n),
\eeq
where $\Pi^+_q(n)$ is the subalgebra of $\Pi_q(n)$ generated by
the variables
$l^+_i[-r]$ with $i=1,\dots,n$ and $r=0,1,\dots$ subject to the relations
$l^+_{i}[0]=1$ for all $i$. Recall the $q$-determinant calculated
in Proposition~\ref{prop:detq}. The following corollary essentially reproduces
the $q$-{\em deformed Miura transformation} of \cite{fr:qa}.

\bco\label{cor:hchy}
We have
\ben
\chi:{\det}_q\big(\Pi+L^+(z)\tss D\tss\de\big)\mapsto
\big(1+\bar\la_1(z)\tss\de\big)\dots \big(1+\bar\la_n(z)\tss\de\big),
\een
where $\bar\la_i(z)=q^{n-2i+1}\ts l^+_i(z)$.
\eco

\bpf
This is immediate from Proposition~\ref{prop:detq} and Theorem~\ref{thm:hch}.
\epf

In the remainder of this section we outline an alternative
construction of the Harish-Chandra homomorphism for
the quantum affine algebra $\U_q(\wh\gl_n)_{\text{\rm cri}}$.
The starting point is the version of the Poincar\'e--Birkhoff--Witt theorem
for a different ordering on the generators as provided by Proposition~\ref{prop:pbwopp}.
The arguments are essentially the same, with only minor changes in notation.
As above, we define the projection
$\theta^{\tss\prime}:\U_q(\wh\gl_n)_{\text{\rm cri}}\to \U^{0}$ in the same way.
Proposition~\ref{prop:hchhom} holds in the same form
but for the different Harish-Chandra homomorphism
\beql{hchprime}
\chi^{\tss\prime}:\Zr_q(\wh\gl_n)\to \wt\Pi_q(n)
\eeq
defined as the restriction of the composition $\chi^{\tss\prime}=\eta\circ\theta^{\tss\prime}$.
For $i=1,\dots,n$ set
\ben
\la'_i(z)=q^{n-2i+1}\ts \frac{l^{+}_{i}(z)\ts l^{-}_{i+1}(zq^{n-2i})\dots l^{-}_{n}(zq^{-n+2})}
{l^{-}_{i}(zq^{n-2i})\dots l^{-}_{n}(zq^{-n})}.
\een
This is a Laurent series in $z$ whose coefficients
are elements of the completed algebra $\wt\Pi_q(n)$.

\bth\label{thm:hchprime}
For each $k=1,\dots,n$ the image of the series $\ell_k(z)$ under the
Harish-Chandra homomorphism \eqref{hchprime} is found by
\ben
\chi^{\tss\prime}:\ell_k(z)\mapsto \sum_{n\geqslant i_1>\dots>i_k\geqslant 1}
\la'_{i_1}(z)\la'_{i_2}(zq^{-2})\dots \la'_{i_k}(zq^{-2k+2}).
\een
\eth

\bpf
The starting point is formula \eqref{lkzhigh} and
the argument is quite similar to the proof of Theorem~\ref{thm:hch}.
\epf

\section{Eigenvalues in Wakimoto modules}
\label{sec:ewm}
\setcounter{equation}{0}

Our goal in this section is to relate the image of the series
$\ell_k(z)$ under the Harish-Chandra homomorphism provided by Theorem~\ref{thm:hch}
with their eigenvalues in the $q$-deformed
Wakimoto modules constructed by Awata, Odake and Shiraishi~\cite{aos:fb}.
Equivalently, due to the work of Frenkel and Reshetikhin~\cite{fr:qa},
these eigenvalues can be interpreted as
elements of the $q$-deformed classical $\Wc$-algebra $\Wc_q(\gl_n)$.
They were associated in \cite{fr:qa} to the series $\ell_{V_k}(z)$
corresponding to fundamental representations
$V_k$ of the quantum affine algebra $\U_q(\wh\sll_n)$.

To establish the relationship, we consider the Wakimoto modules
at the critical level over $\U_q(\wh\gl_n)$.
The coefficients of the series $\ell_k(z)$ act
as multiplications by scalars in the irreducible modules.
We will show that these scalars can be found from Theorem~\ref{thm:hch}
by an appropriate identification of the parameters of the Wakimoto modules
with elements of $\Pi_q(n)$.

The free field realization of \cite{aos:fb} is given in terms of Drinfeld's
``new realization" \cite{d:nr}
of the quantum affine algebra.
Following \cite{fr:qa}, we will use the Ding--Frenkel isomorphism~\cite{df:it}
to get the formulas for the action of the generators
of $\U_q(\wh\gl_n)$ in the Wakimoto modules in terms of the $RLL$ presentation.
Introduce the series $e^{\pm}_{ij}(u)$, $f^{\pm}_{ij}(u)$ and $k^{\pm}_{i}(u)$
which are uniquely determined by the Gauss decompositions of the respective
matrices $L^{\ts\pm}(u)$:
\begin{multline}
L^{\ts\pm}(u)=\\
\begin{bmatrix}1&0&\dots&0\ts\\
                         e^{\pm}_{21}(u)&\ddots&\ddots&\vdots\\
                         \vdots&\ddots&\ddots&0\\
                         e^{\pm}_{n1}(u)&\dots&e^{\pm}_{n\ts n-1}(u)&1
           \end{bmatrix}
           \begin{bmatrix}k^{\pm}_{1}(u)&0&\dots&0\\
                         0&\ddots&\ddots&\vdots\\
                         \vdots&\ddots&\ddots&0\\
                         0&\dots&0&k^{\pm}_{n}(u)
           \end{bmatrix}
           \begin{bmatrix}{\ts 1}&f^{\pm}_{12}(u)&\dots&f^{\pm}_{1n}(u)\\
                         {\ts 0}&\ddots&\ddots&\vdots\\
                         \vdots&\ddots&\ddots&f^{\pm}_{n-1\ts n}(u)\\
                         {\ts 0}&\dots&0&1
           \end{bmatrix}.
\non
\end{multline}
We will need the following quantum minor expressions for these series.
Their Yangian counterparts go back to
\cite{d:nr} and detailed arguments were given in \cite{bk:pp}; see also
\cite[Sec.~1.11]{m:yc}. The quantum
affine algebra case is quite similar so we only sketch the main steps of the proof.

\ble\label{lem:qmin}
We have
\beql{kall}
k^{\pm}_i(u)={L^{\pm}(q^{2i-2}u)\tss}^{1\cdots\ts i}_{1\cdots\ts i}\ts
\Big[{L^{\pm}(q^{2i-2}u)\tss}^{1\cdots\ts i-1}_{1\cdots\ts i-1}\Big]^{-1}
\eeq
for $i=1,\dots,n$ and
\begin{align}
e^{\pm}_{ji}(u)={L^{\pm}(q^{2i-2}u)\tss}^{1\cdots\ts i-1\ts j}_{1\cdots\ts i}\ts
\Big[{L^{\pm}(q^{2i-2}u)\tss}&^{1\cdots\ts i}_{1\cdots\ts i}\Big]^{-1}
\non\\
&{}=\qin\tss \Big[{L^{\pm}(q^{2i}u)\tss}^{1\cdots\ts i}_{1\cdots\ts i}\Big]^{-1}
\ts {L^{\pm}(q^{2i}u)\tss}^{1\cdots\ts i-1\ts j}_{1\cdots\ts i},
\label{ell}
\end{align}
\ben
f^{\pm}_{ij}(u)=\Big[{L^{\pm}(q^{2i-2}u)\tss}^{1\cdots\ts i}_{1\cdots\ts i}\Big]^{-1}
{L^{\pm}(q^{2i-2}u)\tss}^{1\cdots\ts i}_{1\cdots\ts i-1\ts j}
=q\ts {L^{\pm}(q^{2i}u)\tss}^{1\cdots\ts i}_{1\cdots\ts i-1\ts j}
\ts \Big[{L^{\pm}(q^{2i}u)\tss}^{1\cdots\ts i}_{1\cdots\ts i}\Big]^{-1}
\een
for $1\leqslant i<j\leqslant n$.
\ele

\bpf
The arguments for the matrices $L^+(u)$ and $L^-(u)$ are the same so we will
work with $L^+(u)$ and use the notation $\Y_q(\gl_n)$ for the subalgebra
of the quantum affine algebra generated by the coefficients of all series $l^+_{ij}(u)$.
We will also use the algebra $\Y_{\qin}(\gl_n)$ and denote its generator matrix
by $\wt L^{\ts +}(u)=[\ts{\wt l}^{\ts +}_{ij}(u)\tss]$. Due to the property
\ben
R_{21}(v,u)=-R(u,v)\big|_{q\ts\mapsto \qin}
\een
of the $R$-matrix \eqref{ruv}, the mapping
\ben
\om_n: L^+(u)\mapsto \wt L^{\ts +}(u)^{-1}
\een
defines a homomorphism $\om_n:\Y_q(\gl_n)\to \Y_{\qin}(\gl_n)$.
For any $m\geqslant 0$ consider another homomorphism
\beql{phimj}
\jmath_m:\Y_{\qin}(\gl_n)\to \Y_{\qin}(\gl_{m+n})
\eeq
which takes the coefficients of the series
${\wt l}^{\ts +}_{ij}(u)$ to
the respective coefficients of the series
${\wt l}^{\ts +}_{m+i\ts m+j}(u)$.
Consider the composition
$
\phi_m=\omega_{m+n}^{-1}\circ \jmath_m \circ \om_n
$
which is an algebra homomorphism
\ben
\phi_m: \Y_q(\gl_n)\to \Y_q(\gl_{m+n}).
\een
By \cite[Lemma~3.7]{hm:qa}\footnote{The generator matrices $T(u)$ and $\overline T(u)$
of \cite{hm:qa} correspond to $L^-(u)$ and $L^+(u)$, respectively.}
\beql{phim}
\phi_m:l^+_{ij}(u)\mapsto
\Big[{L^{+}(q^{2m}u)\tss}^{1\cdots\ts m}_{1\cdots\ts m}\Big]^{-1}\ts
{L^{+}(q^{2m}u)\tss}^{1\cdots\ts m\ts m+i}_{1\cdots\ts m\ts m+j}.
\eeq
Now apply \cite[Lemmas~1.11.2 and 1.11.5]{m:yc} to the algebra $\Y_q(\gl_n)$
to get
\ben
\bal
k_i^+(u)&=\phi_{i-1}\big(l^+_{11}(u)\big),\\
e^{+}_{ji}(u)&=\phi_{i-1}\big(l^+_{j-i+1\ts 1}(u)\tss l^+_{11}(u)^{-1}\big),\\
f^{+}_{ij}(u)&=\phi_{i-1}\big(l^+_{11}(u)^{-1}\tss l^+_{1\ts j-i+1}(u)\big).
\eal
\een
Together with \eqref{phim} this proves the formula for $k_i^+(u)$ and the first
expressions for $e^{+}_{ji}(u)$ and $f^{+}_{ij}(u)$. To prove the second expressions,
use the relations
\ben
l^+_{j-i+1\ts 1}(u)\tss l^+_{11}(u)^{-1}=\qin\ts l^+_{11}(q^2u)^{-1}\tss l^+_{j-i+1\ts 1}(q^2u)
\een
and
\ben
l^+_{11}(u)^{-1}\tss l^+_{1\ts j-i+1}(u)=q\ts l^+_{1\ts j-i+1}(q^2u)\tss l^+_{11}(q^2u)^{-1}
\een
in $\Y_q(\gl_n)$ implied by \eqref{RLL} (or \eqref{ll}) and apply \eqref{phim} again.
\epf

\bco\label{cor:qdexp}
We have the expansions for the quantum determinants
\ben
\qdet L^{\pm}(z)=k_1^{\pm}(z)\tss k_2^{\pm}(zq^{-2})\dots k_n^{\pm}(zq^{-2n+2}).
\een
\eco

\bpf
This is immediate from Lemma~\ref{lem:qmin}.
\epf

We adopt the notation $\psi_{\pm}^{\tss i}(z)$ and $E^{\pm,i}(z)$ of \cite{aos:fb} for the
series of Drinfeld generators
of $\U_q(\wh\sll_n)$ and identify them with the respective elements of
$\U_q(\wh\gl_n)$ by using the isomorphism of \cite{df:it}
in the normalization of \cite{fm:ha}. The formulas are given in
\cite{fm:ha} only for the zero level case, but they can be easily modified
to include the central element $q^c$ as in \cite{df:it}.
For $i=1,\dots,n-1$ we have
\begin{align}
\psi_{\pm}^{\tss i}(z)&=k^{\mp}_{i+1}(zq^{-i})\ts k^{\mp}_{i}(zq^{-i})^{-1},
\label{psipm}\\[0.4em]
E^{+,i}(z)&=\frac{e^{\ts +}_{i+1\ts i}(zq^{\frac{c}{2}-i})
-e^{\ts -}_{i+1\ts i}(zq^{-\frac{c}{2}-i})}{(q-\qin)\tss z},
\label{epi}\\[0.4em]
E^{-,i}(z)&=\frac{f^{\ts +}_{i\ts i+1}(zq^{-\frac{c}{2}-i})
-f^{\ts -}_{i\ts i+1}(zq^{\frac{c}{2}-i})}{(q-\qin)\tss z}.
\non
\end{align}
In fact, these relations were established in \cite{df:it} for the respective
elements of the algebra $\overline\U_q(\wh\gl_n)$ which we defined in Sec.~\ref{sec:pbw}.
However, the map \eqref{heisext} connecting the two algebras
only results in the multiplication of the series $k^{+}_{i}(z)$
by a scalar series with coefficients in the Heisenberg algebra $\Hc_q(n)$
and in the multiplication of the series $k^{-}_{i}(z)$
by a different scalar series. The series $e^{\pm}_{ij}(u)$ and $f^{\pm}_{ij}(u)$
are not affected and hence the above relations apply to the algebra
$\U_q(\wh\gl_n)$ in the same form.

The $q$-deformed Wakimoto modules over $\U_q(\wh\sll_n)$ are realized in the boson
Fock space by an explicit action of the series $\psi_{\pm}^{\tss i}(z)$ and $E^{\pm,i}(z)$
described in \cite{aos:fb}.
We recall this construction assuming that the level is critical.
That is, we take $k=-n$ in the notation of
\cite{aos:fb} and
for all $1\leqslant i<j\leqslant n$ consider {\em free boson fields}
\ben
\bal
b^{ij}(z)&=-\sum_{r\ne 0}\frac{b_r^{ij}}{[r]}\ts z^{-r}
+\frac{q-\qin}{2\tss\log q}b^{ij}_0\tss\log z+Q^{ij}_b,\\[0.3em]
c^{ij}(z)&=-\sum_{r\ne 0}\frac{c_r^{ij}}{[r]}\ts z^{-r}
+\frac{q-\qin}{2\tss\log q}c^{ij}_0\tss\log z+Q^{ij}_c,
\eal
\een
and
\ben
\bal
b^{ij}_{\pm}(z)&=\pm(q-\qin)\Big(\ts\frac{b^{ij}_0}{2}+\sum_{\pm r> 0}b_r^{ij}\ts z^{-r}\Big),\\[0.3em]
c^{ij}_{\pm}(z)&=\pm(q-\qin)\Big(\ts\frac{c^{ij}_0}{2}+\sum_{\pm r> 0}c_r^{ij}\ts z^{-r}\Big),
\eal
\een
where
\ben
[r]=\frac{q^r-q^{-r}}{q-\qin}
\een
and the coefficients satisfy the relations
\begin{align}
[b_r^{ij},b_s^{kl}]&=-\frac{[r]^2}{r}\ts \de_{ik}\ts\de_{jl}\ts\de_{r,-s},&&
[b_0^{ij},Q_b^{kl}]=-\frac{2\tss\log q}{q-\qin}\ts \de_{ik}\ts\de_{jl},
\label{bq}\\[0.3em]
[c_r^{ij},c_s^{kl}]&=\frac{[r]^2}{r}\ts \de_{ik}\ts\de_{jl}\ts\de_{r,-s},&&
[c_0^{ij},Q_c^{kl}]=\frac{2\tss\log q}{q-\qin}\ts \de_{ik}\ts\de_{jl},
\label{cq}
\end{align}
all other pairs of coefficients commute.
The quantum Heisenberg algebra $\Ac_q(n)$ is generated by the
elements $e^{\pm Q^{ij}_b}$, $e^{\pm Q^{ij}_c}$, $e^{\pm(q-\qin)b^{ij}_0/2}$,
$e^{\pm(q-\qin)c^{ij}_0/2}$ and $b_r^{ij}$, $c_r^{ij}$ with $r\ne 0$.
The defining relations are those implied by the above commutations relations for the
coefficients of the free boson fields.
The Fock representation $F_q(n)$ of $\Ac_q(n)$ is generated by the vacuum vector $|0\rangle$
such that
\ben
b_r^{ij}|0\rangle=c_r^{ij}|0\rangle=0\qquad\text{for all $i<j$ and $r\geqslant 0$}
\een
so that, in particular,
\ben
e^{\pm(q-\qin)b^{ij}_0/2}|0\rangle=e^{\pm(q-\qin)c^{ij}_0/2}|0\rangle=|0\rangle.
\een

The generators $a^i_r$ with $i=1,\dots,n-1$ and $r\in\ZZ$
used in the free field realization of \cite{aos:fb} pairwise commute
at the critical level. Therefore, they may be regarded as
numerical parameters of the $q$-deformed Wakimoto modules and we must have
$a^i_r=0$ for $r\geqslant 0$. For the use in the formulas below we set
\ben
a^i_{\pm}(z)=\pm(q-\qin)\Big(\ts\frac12\ts a^i_0+\sum_{\pm r>0}a^i_r\ts z^{-r}\Big),
\een
as in \cite{aos:fb}, and we have $a^i_{+}(z)=0$ for all $i$.
The ($q$-deformed) Wakimoto module over
$\U_q(\wh\sll_n)$ at the critical level
is defined by the action of the Drinfeld generators
in the space $F_q(n)$. For the action of the coefficients of
the series $\psi_{\pm}^{\tss i}(z)$ we have
\begin{align}
\psi_{\pm}^{\tss i}(z q^{\mp n/2})=
\exp&\Big(\sum_{j=1}^{i}\big(b^{j\ts i+1}_{\pm}(zq^{\pm(j-n-1)})
-b^{j\tss i}_{\pm}(zq^{\pm(j-n)})\big)\non\\
&{}+a^i_{\pm}(z)+\sum_{j=i+1}^{n}\big(b^{i\tss j}_{\pm}(zq^{\pm(j-n)})
-b^{i+1\ts j}_{\pm}(zq^{\pm(j-n-1)})\big)\Big),
\label{kiz}
\end{align}
where $b^{ii}_{\pm}(z)=c^{ii}_{\pm}(z)=0$.
The coefficients of $E^{+,i}(z)$ act by
\begin{align}
E^{+,i}(z)&=\frac{1}{(q-\qin)\tss z}\ts\sum_{j=1}^{i}:\exp\big((b+c)^{j\tss i}(zq^{j-1})\big)
\non\\
&{}\times\Big({\exp}\big(b^{j\ts i+1}_-(zq^{j-1})-(b+c)^{j\ts i+1}(zq^{j-2})\big)
-\exp\big(b^{j\ts i+1}_+(zq^{j-1})-(b+c)^{j\ts i+1}(zq^{j})\big)\Big)
\non\\
&{}\times\exp\Big(\sum_{l=1}^{j-1}\big(b^{l\ts i+1}_+(zq^{l-1})-b^{l\tss i}_+(zq^{l})\big)\Big):,
\label{xizact}
\end{align}
where we have used the notation $(b+c)^{ij}(z)=b^{ij}(z)+c^{ij}(z)$ and set
$b^{ii}(z)=c^{ii}(z)=0$. The colons indicate normal ordering
so that the coefficients $b^{ij}_r$ with $r<0$ or $\exp{Q^{ij}_b}$ should be placed
to the left of the coefficients $b^{ij}_r$ with $r\geqslant 0$. The same rule
applies to the coefficients of $c^{ij}(z)$.
We will not reproduce the formulas for the action
of $E^{-,i}(z)$ as they are given by longer expressions and will not
be used; see \cite[(3.7)]{aos:fb}. Our notation
is the same as in \cite{aos:fb}
for $b_r^{ij}$ and $c_r^{ij}$, whereas $Q^{ij}_b=\hat q^{\tss ij}_b$ and
$Q^{ij}_c=\hat q^{\tss ij}_c$.

By Lemma~\ref{lem:qd}, the coefficients
of the quantum determinants are central in the algebra $\U_q(\wh\gl_n)$
at the critical level. Therefore, the irreducible
Wakimoto modules can be extended to $\U_q(\wh\gl_n)_{\text{\rm cri}}$
by specifying the eigenvalues $K^{\pm}(z)$ of $\qdet L^{\pm}(z)$.
By Corollary~\ref{cor:qdexp}, this gives the conditions
\ben
k_1^{\pm}(z)\tss k_2^{\pm}(zq^{-2})\dots k_n^{\pm}(zq^{-2n+2})\mapsto K^{\pm}(z),
\een
where $K^{+}(z)$ and $K^{-}(z)$ are power series in $z$ and $z^{-1}$, respectively.
Hence, relations \eqref{psipm} and \eqref{kiz} allow us to define the action of the
coefficients of all series $k_i^{\pm}(z)$ on the Fock space.

For any $X\in\U_q(\wh\gl_n)_{\text{\rm cri}}$ we will write $\vacl X\vacr$ to denote
the coefficient of $\vacr$ in the expansion of $X\vacr$ along
the basis of the Fock space.\footnote{An equivalent interpretation would involve
an inner product on the Fock space which we will not introduce.}
More generally, a relation of the form $\vacl X=d\ts\vacl$ for a constant $d$
will be understood in the sense that $\vacl XY\vacr=d\ts\vacl Y\vacr$
for any element $Y\in \U_q(\wh\gl_n)_{\text{\rm cri}}$.
Using this notation
we can parameterize the corresponding modules over $\U_q(\wh\gl_n)_{\text{\rm cri}}$
by the power series
$\vk^+_i(z)$ and $\vk^-_i(z)$ in $z$ and $z^{-1}$, respectively, such that
\beql{kka}
k^{-}_{i}(z)\vacr=\vk^-_i(z)\vacr\Fand \vacl k^{+}_{i}(z)=\vacl\vk^+_i(z)
\eeq
for all $i=1,\dots,n$ satisfying the relations
\ben
\vk^-_{i+1}(z)\vk^-_{i}(z)^{-1}=\exp\big(a^i_{+}(zq^{\frac{n}{2}+i})\big)
\Fand
\vk^+_{i+1}(z)\vk^+_{i}(z)^{-1}=\exp\big(a^i_{-}(zq^{-\frac{n}{2}+i})\big)
\een
for $i=1,\dots,n-1$. Since $a^i_{+}(z)=0$, the series
$\vk^-_{i}(z)$ is the same for each $i$ and
we denote it by $\vk^-(z)$.

The following theorem is essentially due to \cite{fr:qa} subject to the
identification of $\ell_k(z)$ with the series $\ell_{V_k}(z)$
corresponding to the fundamental representation
$V_k$ of $\U_q(\wh\sll_n)$, although the arguments were only outlined there.
The eigenvalues of $\ell_{V_k}(z)$ were interpreted in \cite{fr:qa} as generators
of the $q$-deformed classical $\Wc$-algebras and the Poisson brackets between
the generators were explicitly calculated.

\bth\label{thm:eigen}
Given an irreducible Wakimoto module over $\U_q(\wh\gl_n)_{\text{\rm cri}}$ with
the parameters $\vk^{+}_i(z)$ and $\vk^-(z)$,
the eigenvalues of the coefficients of the series
$\ell_k(z)$ in the module are found by
\ben
\ell_k(z)\mapsto\sum_{1\leqslant i_1<\dots<i_k\leqslant n}
\La_{i_1}(z)\La_{i_2}(zq^{-2})\dots \La_{i_k}(zq^{-2k+2}), \qquad k=1,\dots,n,
\een
where
\ben
\La_i(z)=q^{n-2i+1}\ts \vk^+_{i}(z)\vk^-(zq^{-n})^{-1},\qquad i=1,\dots,n.
\een
\eth

\bpf
Any irreducible Wakimoto module coincides with the cyclic span over $\U_q(\wh\gl_n)_{\text{\rm cri}}$
of the vacuum vector $\vacr$. Hence, the eigenvalues of the coefficients
of the series $\ell_k(z)$ can be found by calculating
the series $\vacl\ell_k(z)\vacr$.

We find from \eqref{xizact} that $E^{+,i}(z)\vacr$ is a power series in $z$ for all $i$.
Therefore, \eqref{epi} implies $e^{\ts -}_{i+1\ts i}(z)\vacr=0$ for $i=1,\dots,n-1$.
Note the following relations for the action of the generators $l^-_{i\tss j}[0]$
on the vacuum vector:
\beql{actvac}
l^-_{i+1\ts i}[0]\vacr =0\Fand l^-_{i\tss i}[0]\vacr = \vk^-[0]\vacr,
\eeq
where $\vk^-[0]$ denotes the constant term of the series $\vk^-(z)$.
Indeed, using \eqref{kall} and \eqref{kka},
we find that ${L^{-}(z)\tss}^{1\cdots\ts i}_{1\cdots\ts i}\vacr$ is a scalar power
series in $z^{-1}$.
Expanding the quantum minor by \eqref{qminorgenl} and taking the constant term
we find that
\beql{qdect}
\sum_{\si\in \Sym_i} (-q)^{-l(\si)} \ts l^{-}_{\si(1)\ts 1}[0]\cdots
l^{-}_{\si(i)\ts i}[0]\vacr
\eeq
is a scalar multiple of the vacuum vector $\vacr$. However, $l^{-}_{j\tss i}[0]=0$
for $j<i$ by \eqref{llt} so that the only nonzero term in \eqref{qdect} corresponds to
the identity permutation $\si$. Therefore, taking the constant terms in
\eqref{kall} and \eqref{kka} we derive the second relation in \eqref{actvac}.
Now use the relation
${L^{-}(z)\tss}^{1\cdots\ts i-1\ts i+1}_{1\cdots\ts i}\vacr=0$
implied by \eqref{ell}.
Exactly as above, we get
\beql{lipp}
l^{-}_{1\tss 1}[0]\cdots l^{-}_{i-1\ts i-1}[0]\ts
l^{-}_{i+1 \ts i}[0]\vacr=0
\eeq
which gives the first relation in \eqref{actvac} by \eqref{lpze}.
As a next step, we will derive the relations
\beql{eji}
e^{\ts -}_{j\ts i}(z)\vacr=0\qquad\text{for all\ \  $j>i$}.
\eeq
It follows from \eqref{RLL} that
\beql{lipone}
\big[l^{-}_{i\ts j}[0],l^{-}_{j\tss m}(z)\big]=(q-\qin)\tss l^{-}_{i\tss m}(z)\ts
l^{-}_{j\tss j}[0]\qquad\text{for \ \ $i>j>m$}
\eeq
and
\ben
\big[l^{-}_{i\ts j}[0],l^{-}_{k\tss m}(z)\big]=0\qquad\text{for \ \ $i,j>k,m$}.
\een
Hence, \eqref{ell} gives
\ben
\big[l^{-}_{j+1\ts j}[0],e^{\ts -}_{j\ts i}(z)\big]=(q-\qin)\tss e^{\ts -}_{j+1\ts i}(z)\ts
l^{-}_{j\tss j}[0]
\een
and \eqref{eji} follows by induction from \eqref{actvac}.

Thus, we may conclude that
\ben
{L^{-}(z)\tss}^{1\cdots\ts i-1\ts j}_{1\cdots\ts i}\vacr=0
\qquad\text{for all\ \ $j>i$}.
\een
Expanding the quantum minor with the use of \eqref{qminorgenl}
and \eqref{qmsym2} we get
\beql{expaji}
\sum_{\si} (-q)^{l(\si)} \ts l^{-}_{\si(j)\tss i}(z)\ts
l^{-}_{\si(i-1)\tss i-1}(q^{-2}z)\cdots
l^{-}_{\si(1)\tss 1}(q^{-2i+2}z)\vacr=0,
\eeq
where the sum is taken over permutations $\si$ of the set $\{1,\dots,i-1,j\}$.
Together with the property that
\beql{expaii}
{L^{-}(z)\tss}^{1\cdots\ts i}_{1\cdots\ts i}\vacr
=\sum_{\si\in\Sym_i} (-q)^{l(\si)} \ts l^{-}_{\si(i)\tss i}(z)\cdots
l^{-}_{\si(1)\tss 1}(q^{-2i+2}z)\vacr
\eeq
is a scalar power series in $z^{-1}$, these relations imply
that
\beql{lmji}
l^{-}_{j\tss i}(z)\vacr=0\qquad\text{for all\ \  $j>i$}
\eeq
and $l^{-}_{i\tss i}(z)\vacr$
is a scalar power series in $z^{-1}$ by an obvious induction. Moreover,
it follows from \eqref{kall} and \eqref{kka} that
\beql{lmii}
l^{-}_{i\tss i}(z)\vacr=\vk^-(z)\vacr\qquad\text{for $i=1,\dots,n$}.
\eeq
Now we will apply similar arguments to derive that
\beql{lva}
\vacl l^{+}_{j\tss i}(z)=0\quad\text{for all\ \  $j>i$}\Fand
\vacl l^{+}_{i\tss i}(z)=\vacl\vk^+_i(z)\quad\text{for\ \  $i=1,\dots,n$}.
\eeq
The first step is to observe that $\vacl z\tss E^{+,i}(z)$
is a power series in $z^{-1}$. Indeed, this follows from
\eqref{xizact} with the use of the relations
\ben
\vacl b_r^{ij}=\vacl c_r^{ij}=0\qquad\text{for all\ \ $r\geqslant 0$}.
\een
One additional step is to use the commutation relations \eqref{bq} and \eqref{cq}
which imply
\ben
\exp Q^{ij}_b
\cdot z^{\frac{q-\qin}{2\tss\log q}\tss b^{ij}_0}
=z^{\frac{q-\qin}{2\tss\log q}\tss b^{ij}_0}\cdot \exp Q^{ij}_b\cdot z
\een
and
\ben
\exp Q^{ij}_c
\cdot z^{\frac{q-\qin}{2\tss\log q}\tss c^{ij}_0}
=z^{\frac{q-\qin}{2\tss\log q}\tss c^{ij}_0}\cdot \exp Q^{ij}_c\cdot z^{-1}.
\een
Although extra powers of $z$ occur as a result of swapping the coefficients,
these powers arising from the coefficients of the series $b^{ij}(z)$ and
$c^{ij}(z)$ cancel each other. Thus, using \eqref{epi}
and noting that the constant term of $e^{+}_{i+1\ts i}(z)$ is zero,
we come to the relation $\vacl e^{+}_{i+1\ts i}(z)=0$ for $i=1,\dots,n-1$.
The rest of the arguments is essentially the same with some obvious adjustments.
In particular, to evaluate the constant term of the power series
$\vacl{L^{+}(z)\tss}^{1\cdots\ts i}_{1\cdots\ts i}$
we use the expansion
\ben
\vacl\sum_{\si\in \Sym_i} (-q)^{l(\si)} l^{+}_{\si(i)\ts i}[0]\cdots\ts l^{+}_{\si(1)\ts 1}[0]
\een
instead of \eqref{qdect} and note that $l^{+}_{j\ts i}[0]=0$ for $j>i$.
Together with the corresponding counterpart of \eqref{lipp} this implies
\ben
\vacl l^+_{i+1\ts i}[0] =0\Fand \vacl l^+_{i\tss i}[0] = \vacl\vk^+_i[0],
\een
where $\vk^+_i[0]$ denotes the constant term of the series $\vk^+_i(z)$.
In the final part we use the relations
\ben
\vacl\sum_{\si} (-q)^{-l(\si)} \ts l^{+}_{\si(1)\tss 1}(z)
\dots l^{+}_{\si(i-1)\tss i-1}(q^{-2i+4}z)\ts l^{+}_{\si(j)\tss i}(q^{-2i+2}z)
=0
\een
with the sum over permutations $\si$ of the set $\{1,\dots,i-1,j\}$, and
\ben
\vacl{L^{+}(z)\tss}^{1\cdots\ts i}_{1\cdots\ts i}
=\vacl\sum_{\si\in\Sym_i} (-q)^{-l(\si)} l^{+}_{\si(1)\tss 1}(z)\ts\dots
l^{+}_{\si(i)\tss i}(q^{-2i+2}z)
\een
instead of \eqref{expaji} and \eqref{expaii}.

Relations \eqref{lmji}, \eqref{lmii}
and \eqref{lva} allow us to conclude that the eigenvalue
$\vacl\ell_k(z)\vacr$ coincides with the image of the series $\ell_k(z)$
under the Harish-Chandra homomorphism calculated in Theorem~\ref{thm:hch}
for the specialization
\ben
l^{+}_{i}(z)=\vk^+_i(z)\Fand l^{-}_{i}(z)=\vk^-(z)
\een
for $i=1,\dots,n$.
Clearly, then $\la_i(z)$ specializes to $\La_i(z)$ and the proof is complete.
\epf

\bre\label{rem:zeromode}
The fact that the eigenvalues of the coefficients of $\ell_k(z)$
in the Wakimoto modules are consistent with the Harish-Chandra images
provided by Theorem~\ref{thm:hch} relies on the properties
\eqref{lmji}, \eqref{lmii} and \eqref{lva}. It was essential for their derivation
that the ``zero mode matrices" $L^+[0]$ and $L^-[0]$ are upper and lower
triangular, respectively. These properties do not hold for the presentation
of the quantum affine algebra used in \cite{fr:qa}, where
the triangularity of the zero mode matrices is opposite.
\qed
\ere


\newpage

\small

\noindent
L.F. \& E.R.:\newline
Laboratoire de Physique Th\'{e}orique LAPTh,
CNRS and Universit\'{e} de Savoie\newline
BP 110, 74941 Annecy-le-Vieux Cedex, France\newline
luc.frappat@lapth.cnrs.fr\newline
eric.ragoucy@lapth.cnrs.fr

\vspace{5 mm}

\noindent
N.J.:\newline
School of Mathematical Sciences\\
South China University of Technology\\
Guangzhou, Guangdong 510640, China\\
and\\
Department of Mathematics\\
North Carolina State University, Raleigh, NC 27695, USA\\
jing@math.ncsu.edu

\vspace{5 mm}

\noindent
A.M.:\newline
School of Mathematics and Statistics\newline
University of Sydney,
NSW 2006, Australia\newline
alexander.molev@sydney.edu.au

\end{document}